\documentclass{amsart}
\usepackage{amssymb}

\usepackage{amscd}

\newtheorem{theorem}{Theorem}[section]
\newtheorem{corollary}[theorem]{Corollary}

\newtheorem{lemma}[theorem]{Lemma}

\theoremstyle{definition}
\newtheorem{definition}{Definition}[section]

\theoremstyle{remark}
\newtheorem{remark}{Remark}[section]

\begin{document}
\title{Bounding the genus of subvarieties of generic hypersurfaces from below}
\author{Herbert Clemens}
\address{Mathematics Department, University of Utah\\
155 So. 1400 East\\
Salt Lake City, UT 84112, USA}
\email{clemens@math.utah.edu}
\date{June, 2002}
\maketitle

\begin{abstract}
A second-order invariant of C. Voisin gives a powerful method for bounding
from below the geometric genus of a k-dimensional subvariety of a degree-d
hypersurface in complex projective n-space. This work uses the Voisin method
to establish a general bound, which lies behind recent results of G.
Pacienza and Z. Ran.
\end{abstract}

\section{Introduction\protect\footnote{%
Partially supported by NSF grant DMS-9970412}}

This paper does nnot seek to improve recently established lower bounds of
the geometric genus of a subvariety of a generic hypersurface in complex
projective space $\Bbb{P}^{n}$, but rather to highlight the fundamental role
in such bounds played by a new second-order invariant introduced by C.
Voisin in \cite{V}. It is Voisin's invariant which has permitted these
recent improvements and the invariant is as interesting in its own right as
are its applications. It is our purpose to here to distinguish its role.

The invariant probes the scheme of lines of high contact with a
hypersurfaces. It is perhaps a bit surprising that lines play any special
role in determining all subvarieties of low genus but it turns out that they
play a central role, a role enhanced by the geometric effect that
second-order variation of the hypersurface has on the lines of high contact.
Roughly the program is as follows. Let $F$ be a generic homogeneous form of
degree $d$ on $\Bbb{P}^{n}$ and let $X_{F}$ denote the corresponding
hypersurface. One wants to bound from below the geometric genus of any $k$
-dimensional smooth variety $Y_{F}$ can be mapped 
\begin{equation*}
f:Y_{F}\rightarrow X_{F}
\end{equation*}
generically. Such a bound will of course depend on $n$, $d$, and $k$. A
rather straightforward argument by adjunction is available when $d$ is
fairly large with respect to $n$ but the situation is more delicate for
smaller $d$. As $d$ decreases, adjunction is not enough to produce global
section of the canonical bundle of $Y_{F}$. Using this very weakness to
advantage, at a general point $y\in Y_{F}$ Voisin uses the ``possible
degeneracy of adjunction'' to canonically distinguish a line $l\left(
y\right) $among those passing through a point
\begin{equation*}
x=x\left( y\right) \in X_{F}\text{.}
\end{equation*}
She then builds a second-order invariant out of the Lie bracket restricted
to a distinguished distribution inside $T_{Y},$ a kind of  ``second
fundamental form'' with values in $d$-forms on $l\left( y\right) $. Using
this tool, one shows in many situations, that if $Y_{F}$ has low (e.g. zero)
genus, $l\left( y\right) $ must have contact order at least $d$ with $X_{F}$
at $x\left( y\right) $. Adjunction then imples, in the range of values
considered, that such a $Y_{F}$ must actually lie inside the locus of those
lines which lie entirely inside $X_{F}$.

A rather mysterious point in the development of this program is how lines
enter the picture in the first place. The idea is the following. A subbundle
of $T_{X}$ over which one has some control is the so-called ``vertical
tangent space,'' that is, the vectors at a point $\left( x,F\right) \in X$
corresponding to directions in which $F$ moves but $x$ stays fixed. This is
just the space $M_{x}^{d}$ of homogeneous forms $G$ vanishing at $x$. Thus
the critical issue is the positivity of the bundle $M_{\Bbb{P}^{n}}^{d}$
whose fiber at $x$ is $M_{x}^{d}$. Since one is studying adjunction relating
the canonical bundle of $X_{F}$ to that of $Y_{F}$, one studies, for general 
$y\in Y_{F}$ and $x=x\left( y\right) $, the mapping 
\begin{equation}
\mu :M_{x}^{d}\rightarrow \frac{\left. T_{X}\right| _{x}}{f_{*}\left( \left.
T_{Y}\right| _{y}\right) }  \label{zero}
\end{equation}
where $Y$ is the versal family of deformations of $Y_{F}$. A sub-bundle of $%
M_{\Bbb{P}^{n}}^{d}$ whose positivity (and hence global generation) is
easily controlled is 
\begin{equation*}
P_{1}\cdot M_{\Bbb{P}^{n}}^{1}
\end{equation*}
where $P_{1}$ is a (generically chosen) fixed homogeneous form of degree $d-1
$. However the map 
\begin{equation*}
M_{x}^{1}\overset{\cdot P_{1}}{\longrightarrow }\frac{\left. T_{X}\right|
_{x}}{f_{*}\left( \left. T_{Y}\right| _{y}\right) }
\end{equation*}
is not surjective necessarily. If a second map 
\begin{equation*}
M_{x}^{1}\overset{\cdot P_{2}}{\longrightarrow }\frac{\left. T_{X}\right|
_{x}}{f_{*}\left( \left. T_{Y}\right| _{y}\right) +\mu \left( P_{1}\cdot %
M_{x}^{1}\right) }
\end{equation*}
is required to achive surjectivity, the positivity or global-generation
conclusion that can be drawn is weaker. In fact, the large the number of
linearly independent $P$'s one needs to the filling up of the image of $%
\left( \ref{zero}\right) $, the weaker the bounds one gets. On the other
hand, for each additional $P_{s}$, the ranks of the successive maps 
\begin{equation*}
M_{x}^{1}\rightarrow \frac{\left. T_{X}\right| _{x}}{f_{*}\left( \left.
T_{Y}\right| _{y}\right) +\mu \left( P_{1}\cdot M_{x}^{1}+\ldots +P_{s-1}%
\cdot M_{x}^{1}\right) }
\end{equation*}
are non-increasing. Lines enter for $n$, $d$, $k$ for which $s$ is large
enough that the last two maps in this sequence have rank $1$, which turns
out to be exactly the range is which positivity is insufficient to apply
adjunction directly.. In this case one associates to $x=x\left( y\right) $
the line $l\left( y\right) $ whose ideal is given by the kernel of the
(first) rank-$1$ mapping.

The appearance of the rational mapping 
\begin{equation*}
l:Y\rightarrow G
\end{equation*}
where $G$ is the Grassmann of projective lines in $\Bbb{P}^{n}$ allows us to
lift the mapping $f$ to a (generically injective) mapping 
\begin{equation*}
g=\left( l,x,F\right) :Y\rightarrow \Delta :=L\times _{\Bbb{P}^{n}}X
\end{equation*}
where
\begin{equation}
\begin{array}{ccc}
L & \overset{q}{\longrightarrow } & \Bbb{P}^{n} \\ 
\downarrow ^{p} &  &  \\ 
G &  & 
\end{array}
\label{univ}
\end{equation}
is the universal line. This line in turn determines a distinguished
(vertical) distribution in $T_{Y}.$ To describe it, let $G_{x}\subseteq G$
denote the variety of lines through $x$. For $l\in G$, let $M_{l}^{d}$
denote homogeneous forms of degree $d$ which vanish on $l$. Then one can
write 
\begin{equation*}
\left. T_{\Delta }\right| _{\left( l,x,F\right) }=\left. T_{G_{x}}\right|
_{l}\oplus \left. T_{\Bbb{P}^{n}}\right| _{x}\oplus M_{x}^{d}
\end{equation*}
and consider the subspace 
\begin{equation*}
\left. T^{\prime }\right| _{y}:=\left( 0\oplus 0\oplus M_{l}^{d}\right) \cap
g_{*}\left( \left. T_{Y}\right| _{y}\right) 
\end{equation*}
of 
\begin{equation*}
\left. T_{\Delta }\right| _{\left( l\left( y\right) ,x\left( y\right)
,F\left( y\right) \right) }.
\end{equation*}
The key point which requires Voisin's second-order construction is the proof
that the distribution $T^{\prime }$ is integrable. Said another way, one
shows that, when $y$ moves (infinitesimally) so that the derivative of $%
F\left( y\right) $ is a polynomial vanishing on $l\left( y\right) $, then,
in that direction, the line $l\left( y\right) $is (infinitesimally)
stationary. This, the fact that $Y$ can be taken to be $GL\left( n+1\right) $%
-invariant, and a little algebra yields the conclusion that $l\left(
y\right) $ must have contact at least $d$ with $X_{F\left( y\right) }$ at $%
x\left( y\right) $.

The construction of the second-order invariant which shows the integrability
of the distribution $T^{\prime }$ goes roughly as follows. Referring to $%
\left( \ref{univ}\right) $, let $E^{d}$ denote the bundle 
\begin{equation*}
p_{*}q^{*}\mathcal{O}_{\Bbb{P}^{n}}\left( d\right) 
\end{equation*}
on $G$ and let
\begin{equation*}
E_{\Bbb{P}^{n}}^{d}\subseteq p^{*}E^{d}
\end{equation*}
denote the sub-bundle whose fiber at $\left( l,x\right) $ is given by the
set of $d$-forms on $l$ which vanish at $x.$ Now at a general point
\begin{equation*}
\left( l,x,F\right) =\left( l\left( y\right) ,x\left( y\right) ,F\left(
y\right) \right) 
\end{equation*}
in the image of $Y$ in $\Delta $, $M_{x}^{d}$ sits as a summand of $\left.
T_{\Delta }\right| _{\left( l,x,F\right) }$ and maps to the fiber $E_{\left(
l,x\right) }^{d}$ of $E_{\Bbb{P}^{n}}^{d}$ by evaluation. Voisin constructs
her second-order invariant out of the composition
\begin{equation}
\left[ T^{\prime },T^{\prime }\right] \rightarrow E_{\Bbb{P}^{n}}^{d}
\label{Vmap}
\end{equation}
of Lie bracket in $g_{*}T_{Y}$ with this evalution map, giving it a
beautiful geometric interpretation. This allows the conclusion that the
co-rank of $\left( \ref{Vmap}\right) $ is incompatible with previously
established bounds unless $\left( \ref{Vmap}\right) $ is zero, which in turn
implies that the distribution $T^{\prime }$ is integrable.

G. Pacienza \cite{P} and Z. Ran \cite{R} have taken the Voisin technique
even further and independently produce stronger results than those given in
Main Theorem and Corollary below for subvarieties $Y_{F}$ such that, for
some $a\geq -1$, 
\begin{equation*}
h^{0}\left( \omega _{Y_{F}}\left( a\right) \right) =0.
\end{equation*}
Though earlier versions of this manuscript contributed to the stronger
results of both authors, its primary justification at this point is its
focus on Voisin's new second-order invariant..

\subsection{The formal setting}

Before we can state the main theorem, we need to complete the notation we
will need throughout the paper. For fixed $d$ and $n$ and variable $r$ let 
\begin{eqnarray*}
M^{r} &=&H^{0}\left( \mathcal{O}_{\Bbb{P}^{n}}\left( r\right) \right)  \\
S &=&H^{0}\left( \mathcal{O}_{\Bbb{P}^{n}}\left( d\right) \right) -\left\{
0\right\} .
\end{eqnarray*}
We let 
\begin{equation*}
N=\left( 
\begin{array}{c}
n+d \\ 
d
\end{array}
\right) =\dim S.
\end{equation*}
Again as above let 
\begin{equation*}
X\subseteq \Bbb{P}^{n}\times S
\end{equation*}
denote the universal hypersurface and 
\begin{equation*}
\Delta =L\times _{\Bbb{P}^{n}}X\subseteq G\times \Bbb{P}^{n}\times S
\end{equation*}
the incidence variety. We have the commutative diagram of natural maps: 
\begin{equation}
\begin{array}{ccccccccc}
&  &  &  & \Delta  &  &  &  &  \\ 
&  &  & \swarrow _{\pi } &  & \searrow ^{\rho } &  &  &  \\ 
&  & L &  &  &  & X &  &  \\ 
& \swarrow _{p} &  & \searrow ^{q} &  & \swarrow _{t} &  & \searrow ^{s} & 
\\ 
G &  &  &  & \Bbb{P}^{n} &  &  &  & S
\end{array}
\label{bigdiagram}
\end{equation}
As above we put
\begin{eqnarray*}
M_{\Bbb{P}^{n}}^{r} &=&\left\{ \left( x,P\right) \in \Bbb{P}^{n}\times 
M^{r}:P\left( x\right) =0\right\}  \\
M_{x}^{r} &=&\left. M_{\Bbb{P}^{n}}^{r}\right| _{x} \\
M_{G}^{r} &=&\left\{ \left( l,P\right) \in G\times M^{r}:\left. P\right|
_{l}=0\right\}  \\
M_{l}^{r} &=&\left. M_{G}^{r}\right| _{l}.
\end{eqnarray*}
Then for 
\begin{equation*}
E^{r}=p_{*}\circ q^{*}\mathcal{O}_{\Bbb{P}^{n}}\left( r\right) 
\end{equation*}
we have the exact sequence 
\begin{equation}
0\rightarrow M_{G}^{r}\rightarrow G\times M^{r}\overset{e^{r}}{%
\longrightarrow }E^{r}\rightarrow 0.  \label{es}
\end{equation}

\subsection{Versal subvarieties}

\begin{definition}
Call we call $Y/S^{\prime }$ a \textit{versal sub-family of} $k$\textit{%
-folds} \textit{of }$X/S$ if
\begin{equation*}
Y\subseteq \Bbb{P}^{n}\times S^{\prime }
\end{equation*}
is a smooth and projective family over $S^{\prime }$ of fiber dimension $k$
admitting a $GL\left( n+1\right) $-action such that, for some etale map
\begin{equation*}
S^{\prime }\rightarrow S,
\end{equation*}
there is a generically injective, $GL\left( n+1\right) $-equivariant map 
\begin{equation*}
f:Y/S^{\prime }\rightarrow X/S.
\end{equation*}
\end{definition}

Here we use the word ``versal'' to suggest the fact that, by assumption, $%
Y_{F}$ deforms with every local deformation of $F$. Also the consideration
of $S^{\prime }/S$ etale is necessitated by Stein factorization. Several of
the subvarieties we are studying may occur on the same (general) $X_{F}$,
for example, lines on the general cubic surface. In order that the notation
not become even more combersome we shall to refer to any of these components
as $Y_{F}$. 

For any integer $a\geq 0$ let 
\begin{eqnarray*}
\omega _{X}\left( a\right)  &=&\omega _{X}\otimes t^{*}\mathcal{O}_{\Bbb{P}%
^{n}}\left( a\right)  \\
\omega _{Y}\left( a\right)  &=&\omega _{Y}\otimes \left( t\circ f\right) ^{*}%
\mathcal{O}_{\Bbb{P}^{n}}\left( a\right) .
\end{eqnarray*}
We will be interested in the case in which 
\begin{equation}
h^{0}\left( \omega _{Y_{F}}\left( a\right) \right) =0  \label{fund}
\end{equation}
where 
\begin{equation*}
Y_{F}
\end{equation*}
is a fiber of $Y/S^{\prime }$ lying over a generic $F\in S$.

\subsection{The Main Theorem}

The purpose of this paper is to prove:

\begin{theorem}
\label{MT}Suppose $Y/S^{\prime }$ is a versal subfamily of $X/S$ of $k$
-folds such that $\left( \ref{fund}\right) $ holds for some $a\geq 0$. If
the inequalities 
\begin{equation*}
d+a\geq \max \left\{ \frac{7n-3k-3}{4},\frac{3n-k+1}{2}\right\} 
\end{equation*}
and 
\begin{equation*}
\frac{d\left( d+1\right) }{2}\geq 3n-k-1
\end{equation*}
are satisfied, then the image of $Y/S^{\prime }$ lies inside the locus of
lines on $X/S$.
\end{theorem}

Notice that, if the codimension of $Y$ in $X$ is at most $4$ then the first
inequality in the theorem becomes 
\begin{equation*}
d+a\geq \frac{3n-k+1}{2}
\end{equation*}
while, if the codimension is at least $4$, it becomes 
\begin{equation*}
d+a\geq \frac{7n-3k-3}{4}.
\end{equation*}

\subsection{Corollary of the proof of the Main Theorem}

Let 
\begin{equation*}
\iota _{F}=\left( t\circ \left. f\right| _{Y_{F}}\right) :Y_{F}\rightarrow 
\Bbb{P}^{n}.
\end{equation*}
Since $S^{\prime }/S$ is etale, the natural mapping 
\begin{equation}
\left. T_{S}\right| _{F}\cong \left. T_{S^{\prime }}\right| _{F}\rightarrow
H^{0}\left( N_{\iota _{F}}\right)   \label{localY}
\end{equation}
composes with
\begin{equation*}
H^{0}\left( N_{\iota _{F}}\right) \otimes \mathcal{O}_{Y_{F}}\rightarrow
N_{\iota _{F}}
\end{equation*}
to give a map 
\begin{equation}
\nu :M^{d}\otimes \mathcal{O}_{Y_{F}}\rightarrow N_{\iota _{F}}
\label{numap}
\end{equation}
which is generically surjective by $GL\left( n+1\right) $-equivariance.
Considering $M^{d}$ as a subspace of $\left. T_{\Bbb{P}^{n}\times S}\right|
_{\left( x,F\right) }$, we have for generic $x\in f\left( Y_{F}\right) $
that 
\begin{equation}
F^{\prime }\in M^{d}\cap \left. T_{Y}\right| _{x}\Leftrightarrow \left. \nu
\left( F^{\prime }\right) \right| _{x}=0.  \label{cond}
\end{equation}
An analysis of the proof of Theorem \ref{MT} yields the following.

\begin{corollary}
\label{MC}Suppose that $d\left( d+1\right) \geq 4n+4$ and $\left( \ref{fund}%
\right) $ holds for some $a\geq 0$. Let $Y/S^{\prime }$ be a versal
subfamily of $X/S$ of $\left( n-3\right) $-folds which does not lie inside
the locus of lines on $X/S.$ Then, for generic $P\in M^{d-1}$, the map 
\begin{equation*}
\nu ^{\prime }=\left. \nu \right| _{P\cdot M^{1}}:M^{1}\otimes \mathcal{O}%
_{Y_{F}}\rightarrow N_{\iota _{F}}
\end{equation*}
is generically surjective and so $\nu ^{\prime }$ induces a rational map 
\begin{equation*}
\mu :Y_{F}\rightarrow Gr\left( n-2,n+1\right) .
\end{equation*}
Furthermore 
\begin{equation*}
h^{0}\left( \omega _{Y_{F}}\left( n+2-d\right) \right) \geq \left( 
\begin{array}{c}
n+1 \\ 
n-2
\end{array}
\right) -r
\end{equation*}
where 
\begin{equation}
r=\dim \left\{ h\in H^{0}\left( \mathcal{O}_{Gr\left( n-2,n+1\right) }\left(
1\right) \right) :\left. h\right| _{\mu \left( Y_{K}\right) }=0\right\} .
\label{rdef}
\end{equation}
\end{corollary}

For example, if $X_{F}$ is a sextic threefold and $k=1$, then 
\begin{equation*}
h^{0}\left( \omega _{Y_{F}}\right) =10-r
\end{equation*}
and so 
\begin{equation*}
h^{0}\left( \omega _{Y_{F}}\right) >1
\end{equation*}
unless $\mu \left( Y_{F}\right) $ is a point.

\subsection{Notation for relative tangent spaces}

An important piece of notation throughout this paper is that for the
``vertical tangent space'' $T_{g}$ to a smooth (surjective) map 
\begin{equation*}
g:Y\rightarrow Z.
\end{equation*}
$T_{g}$ is defined by the exact sequence 
\begin{equation*}
0\rightarrow T_{g}\rightarrow T_{Y}\rightarrow g^{*}T_{Y}\rightarrow 0.
\end{equation*}
For a composition 
\begin{equation*}
W\overset{f}{\longrightarrow }Y\overset{g}{\longrightarrow }Z
\end{equation*}
where $g$ and $g\circ f$ are both smooth (surjective), the Snake Lemma gives
the exact sequence 
\begin{equation*}
0\rightarrow T_{f}\rightarrow T_{g\circ f}\rightarrow f^{*}T_{g}\rightarrow
0.
\end{equation*}
Also, for a fibered product of smooth morphisms 
\begin{equation*}
\begin{array}{ccccc}
&  & X\times _{Z}Y &  &  \\ 
& \swarrow _{a} &  & \searrow ^{a^{\prime }} &  \\ 
X &  &  &  & Y \\ 
& \searrow ^{b} &  & \swarrow _{b^{\prime }} &  \\ 
&  & Z &  & 
\end{array}
\end{equation*}
we have isomorphisms 
\begin{eqnarray*}
T_{a} &\rightarrow &T_{b\circ a}\rightarrow \left( a^{\prime }\right)
^{*}T_{b^{\prime }} \\
T_{a^{\prime }} &\rightarrow &T_{b^{\prime }\circ a^{\prime }}\rightarrow
a^{*}T_{b}
\end{eqnarray*}
so that 
\begin{equation*}
T_{b\circ a}=a^{*}T_{b}\oplus \left( a^{\prime }\right) ^{*}T_{b^{\prime }}.
\end{equation*}
For example referring to $\left( \ref{bigdiagram}\right) $%
\begin{equation}
T_{q\circ \pi }=\pi ^{*}T_{q}\oplus \left( t\circ \rho \right) ^{*}M_{\Bbb{P}%
^{n}}^{d}.  \label{sum}
\end{equation}

\section{Conditions for lifting a versal family $Y/S^{\prime }$}

\subsection{Positivity results\label{positivity}}

\begin{lemma}
\label{LP}i) The sheaf 
\begin{equation*}
\Omega _{\Bbb{P}^{n}}^{h}\left( h+1\right) 
\end{equation*}
is generated by global sections.

ii) 
\begin{equation*}
M_{\Bbb{P}^{n}}^{1}=\Omega _{\Bbb{P}^{n}}^{1}\left( 1\right) 
\end{equation*}

so the sheaf 
\begin{equation*}
\left( \bigwedge\nolimits^{h}M_{\Bbb{P}^{n}}^{1}\right) \left( 1\right) 
\end{equation*}
is generated by global sections for all $h\geq 0$.

iii) 
\begin{equation*}
M_{\Bbb{P}^{n}}^{r}\otimes \mathcal{O}_{\Bbb{P}^{n}}\left( 1\right) 
\end{equation*}
is generated by global sections for all $r\geq 0$.

iv) The sheaf 
\begin{equation*}
M_{G}^{r}\left( 1\right) 
\end{equation*}
is generated by global sections for all $r\geq 0$.
\end{lemma}

\begin{proof}
i) This is a corollary of Mumford's $m$-regularity theorem as follows. By
the classical vanishing results of Bott (see \cite{B}, page 246), 
\begin{equation*}
H^{i}\left( \Omega _{\Bbb{P}^{n}}^{h}\left( k^{\prime }\right) \right) =0
\end{equation*}
unless 
\begin{eqnarray*}
i &=&h,\;k^{\prime }=0 \\
i &=&0,\;k^{\prime }>h \\
i &=&n,\;k^{\prime }<h-n
\end{eqnarray*}
so, in particular, for $i>0$ and $k^{\prime }\geq h+1-i$. Therefore, by
Mumford's regularity theorem (\cite{M}, page 99), the maps 
\begin{equation*}
H^{0}\left( \Omega _{\Bbb{P}^{n}}^{h}\left( k^{\prime }\right) \right)
\otimes H^{0}\left( \mathcal{O}_{\Bbb{P}^{n}}\left( 1\right) \right)
\rightarrow H^{0}\left( \Omega _{\Bbb{P}^{n}}^{h}\left( k^{\prime }+1\right)
\right)
\end{equation*}
are surjective for $k^{\prime }\geq h+1$.

ii) The isomorphism 
\begin{equation*}
M_{\Bbb{P}^{n}}^{1}=\Omega _{\Bbb{P}^{n}}^{1}\left( 1\right) .
\end{equation*}
is immediate from the (dual of the) Euler sequence for the tangent bundle of
projective space. Then use i).

iii) Use the surjection 
\begin{equation*}
M^{r-1}\otimes M_{\Bbb{P}^{n}}^{1}\rightarrow M_{\Bbb{P}^{n}}^{r}.
\end{equation*}

iii) Again it suffices to check the case $r=1$. Using the irreducibility of
the action of $GL\left( n+1\right) $ it suffices to construct a single
non-trivial meromorphic section of $M_{G}^{1}$ with simple pole along the
zero set of a Pl\"{u}cker coordinate. To do this, for all lines $l$ not
meeting 
\begin{equation*}
X_{1}=X_{2}=0
\end{equation*}
the Pl\"{u}cker coordinate $p_{12}\left( l\right) \neq 0$ so there is, by
Cramer's rule, a unique $A=X_{0}+aX_{1}+bX_{2}\in M_{G}^{1}$ containing $l$.
\end{proof}

\subsection{Global generation in the vertical tangent space}

Next we wish to study how close 
\begin{equation*}
\bigwedge\nolimits^{h}M_{\Bbb{P}^{n}}^{d}
\end{equation*}
is to being generated by global sections. Fix $x\in \Bbb{P}^{n}$. Let 
\begin{equation*}
T\subseteq M_{x}^{r}
\end{equation*}
be a subspace of corank $h$. For each $s\leq h,$ consider the map 
\begin{equation}
\mu _{T,s}:\sum\nolimits_{m=1}^{s}t^{*}M_{\Bbb{P}^{n}}^{1}\cdot
P_{m}\rightarrow \frac{M_{x}^{r}}{T}.  \label{mumap}
\end{equation}
for generically chosen $P_{1},\ldots ,P_{s}\in M^{r-1}$. This map is
surjective for $s=h$ and its rank 
\begin{equation*}
\gamma \left( s\right)
\end{equation*}
must be a strictly increasing function of $s$ until surjectivity is reached.
Furthermore 
\begin{equation}
\gamma \left( s\right) -\gamma \left( s-1\right)  \label{gammamap}
\end{equation}
is a non-increasing function. Let $s_{T}\leq h$ denote the smallest $s$ such
that $\mu _{T,s}$ is surjective and $s_{T}^{\prime }\leq s_{T}$ be the
smallest value such that 
\begin{equation*}
\gamma \left( s+1\right) -\gamma \left( s\right) \leq 1.
\end{equation*}
Either

i) $s_{T}^{\prime }=s_{T}$ and 
\begin{equation*}
s_{T}\leq \frac{h}{2}
\end{equation*}
or

ii). $s_{T}^{\prime }<s_{T}$ and for generic $P\in M^{r-1}$%
\begin{equation}
\nu _{T,P}:M_{\Bbb{P}^{n}}^{1}\overset{\cdot P}{\longrightarrow }\frac{
t^{*}M_{\Bbb{P}^{n}}^{r}}{T+image\left( \mu _{T,s_{T}^{\prime }}\right) }
\label{map'}
\end{equation}
has rank $1$. Also 
\begin{equation}
h-\left( s_{T}-s_{T}^{\prime }\right) =\dim \left( image\left( \mu
_{T,s_{T}^{\prime }}\right) \right) .  \label{ineq}
\end{equation}

\begin{remark}
\label{list}Notice that, by this last equality, 
\begin{equation*}
\left( s_{T}-s_{T}^{\prime }\right) \geq rank\left( \mu _{T,s_{T}^{\prime
}}\right) 
\end{equation*}
if and only if 
\begin{equation*}
\frac{h}{2}\leq \left( s_{T}-s_{T}^{\prime }\right) .
\end{equation*}
Of course we always have 
\begin{equation*}
s_{T}+s_{T}^{\prime }=\left( s_{T}-s_{T}^{\prime }\right) +2s_{T}^{\prime }%
\leq h.
\end{equation*}
\end{remark}

\begin{lemma}
\label{goodlem}i) If 
\begin{equation*}
s_{T}^{\prime }+1<s_{T},
\end{equation*}
there exists $l_{T}\in G$ such that 
\begin{equation*}
M_{l_{T}}^{r}\subseteq T+image\left( \mu _{T,s_{T}^{\prime }}\right) .
\end{equation*}

ii) The line $l_{T}$ in i) is independent of the choice of $P_{1},\ldots
,P_{s_{T}^{\prime }}$.

iii) More generally (even when $s_{T}^{\prime }\geq s_{T}-1$), given any $l$
such that 
\begin{equation*}
M_{l_{T}}^{r}\subseteq T+image\left( \mu _{T,s_{T}^{\prime }}\right) ,
\end{equation*}
we have 
\begin{equation*}
\dim \frac{T}{T\cap M_{l_{T}}^{r}}=\dim \frac{T+M_{l_{T}}^{r}}{M_{l_{T}}^{r}}%
=r-s_{T},
\end{equation*}
and 
\begin{equation*}
\dim \frac{M_{l_{T}}^{r}}{T\cap M_{l_{T}}^{r}}=\dim \frac{T+M_{l_{T}}^{r}}{T}%
=rank\left( \mu _{T,s_{T}^{\prime }}\right) -s_{T}^{\prime }.
\end{equation*}
\end{lemma}

\begin{proof}
i) The rank of 
\begin{equation*}
\nu _{P_{s}}:M_{\Bbb{P}^{n}}^{1}\overset{\cdot P}{\longrightarrow }\frac{
t^{*}M_{\Bbb{P}^{n}}^{r}}{T+image\left( \mu _{T,s-1}\right) }
\end{equation*}
is locally constant as we vary generic $P_{s}\in M^{r-1}$. Thus 
\begin{equation*}
\frac{\partial \nu _{P_{s}}}{\partial P_{s}}\in Hom\left( kernel\left( \nu
_{P_{s}}\right) ,\ image\left( \nu _{P_{s}}\right) \right)
\end{equation*}
so that 
\begin{equation*}
Q\cdot K_{s}\subseteq image\left( \nu _{P_{s}}\right)
\end{equation*}
for each $Q\in M^{r-1}$. So 
\begin{equation}
\frac{K_{s}\cdot M^{r-1}+T+image\left( \mu _{T,s-1}\right) }{T+image\left(
\mu _{T,s-1}\right) }\subseteq image\left( \nu _{P_{s}}\right) .
\label{star}
\end{equation}
Leaving $x$ and $T$ fixed but varying $P_{s}\in M^{r-1}$, suppose that $%
K_{s} $ varies. Then differentiating $\left( \ref{star}\right) $ with
respect to $P_{s}$, constant rank again implies that 
\begin{equation}
\frac{\left( \frac{\partial K_{P_{s}}}{\partial P_{s}}\right) \cdot
M^{r-1}+T+image\left( \mu _{T,s-1}\right) }{T+image\left( \mu
_{T,s-1}\right) }\subseteq image\left( \nu _{P_{s}}\right) .  \label{star'}
\end{equation}

Now use that 
\begin{equation*}
rank\left( \mu _{T,s_{T}^{\prime }+1}\right) =1
\end{equation*}
so that 
\begin{equation*}
\left. K_{P_{s_{T}^{\prime }+1}}\right| _{x}=M_{l_{T}}^{1}
\end{equation*}
for some line $l_{T}$ passing through $x$. To see that $l$ does not depend
on the (generic) choice of $P_{s_{T}^{\prime }+1}$, notice that $\left( \ref
{star'}\right) $ and $s_{T}^{\prime }+1<s_{T}$ imply that 
\begin{equation*}
M_{l_{T}}^{1}\subseteq \left. K_{P_{s_{T}^{\prime }+1}}+\frac{\partial
K_{P_{s_{T}^{\prime }+1}}}{\partial P_{s}}\right| _{x}\neq M_{x}^{1}
\end{equation*}
so that 
\begin{equation*}
\frac{\partial K_{P_{s_{T}^{\prime }+1}}}{\partial P_{s}}\subseteq
M_{l_{T}}^{1}.
\end{equation*}
So 
\begin{equation*}
\frac{M_{l_{T}}^{1}\cdot M^{r-1}+T}{T}\subseteq
\bigcap\nolimits_{P_{s_{T}^{\prime }+1}}image\left( \mu _{T,s_{T}^{\prime
}+1}\right) =image\left( \mu _{T,s_{T}^{\prime }}\right)
\end{equation*}
since $image\left( \mu _{T,s_{T}^{\prime }+1}\right) $ cannot be independent
of $P_{s_{T}^{\prime }+1}$ and 
\begin{equation*}
\dim \frac{image\left( \mu _{T,s_{T}^{\prime }+1}\right) }{image\left( \mu
_{T,s_{T}^{\prime }}\right) }=1.
\end{equation*}
Thus 
\begin{equation*}
\frac{M_{l_{T}}^{r}+T}{T}\subseteq image\left( \mu _{T,s_{T}^{\prime
}}\right) .
\end{equation*}

For ii), notice that $l_{T}$ is the \textit{unique} line $l$ such that 
\begin{equation*}
M_{l}^{1}\cdot M^{r-1}\subseteq T+M^{1}\cdot P_{1}+\ldots +M^{1}\cdot
P_{s_{T}^{\prime }}+M^{1}\cdot P_{s_{T}^{\prime }+1}
\end{equation*}
and so, for $s\leq s_{T}^{\prime }$, we can reverse the roles of $P_{s}$ and 
$P_{s_{T}^{\prime }+1}$ in $\left( \ref{map'}\right) $ and we must get the 
\textit{same} line $l$ defining the kernel. By repeating this argument for
each $s\leq s_{T}^{\prime }$ we see that the line $l_{T}$ is also
independent of the choice of $P_{1},\ldots ,P_{s_{T}^{\prime }}$.

For iii), let $l$ be any line such that 
\begin{equation*}
\frac{M_{l}^{r}+T}{T}\subseteq image\left( \mu _{T,s_{T}^{\prime }}\right) .
\end{equation*}
Let 
\begin{equation*}
E_{l}^{r}=\frac{M_{x}^{r}}{M_{l}^{r}}
\end{equation*}
and denote the image of $T$ in $E_{l}^{r}$ as $\hat{T}$. Now notice that,
for $P_{s}$ generic, the image $M_{x}^{1}\cdot P_{s}$ in 
\begin{equation}
\frac{E_{l}^{r}}{\hat{T}}  \label{hatquot}
\end{equation}
has rank $1$ for each $s$ until we reach the situation in which 
\begin{equation*}
M_{l}^{r}+image\left( \mu _{s}\right)
\end{equation*}
generates $\left( \ref{hatquot}\right) $. Since 
\begin{equation*}
M_{l}^{r}\subseteq image\left( \mu _{s_{T}^{\prime }}\right)
\end{equation*}
and $s_{T}^{\prime }\leq s_{T}$, this cannot happen until $s=s_{T}$. So 
\begin{equation*}
s_{T}=\dim \frac{M_{x}^{r}}{T+M_{l}^{r}}.
\end{equation*}

The second equality follows since 
\begin{equation*}
\frac{image\left( \mu _{s}\right) }{image\left( \mu _{s-1}\right) }
\end{equation*}
must have a codimension $1$ subspace generated by 
\begin{equation*}
\frac{M_{l}^{r}+T}{T}\cap image\left( \mu _{s}\right)
\end{equation*}
for each $s\leq s_{T}^{\prime }$.
\end{proof}

\subsection{Geometry of versal families}

We are now ready to apply the results of \S \ref{positivity} in the case of
a versal family $Y/S^{\prime }$ with 
\begin{eqnarray*}
r &=&d \\
h &=&\left( n-1\right) -k \\
T &=&f_{*}T_{t\circ f}\subseteq \left. t^{*}M_{\Bbb{P}^{n}}^{d}\right|
_{\left( x,F\right) }
\end{eqnarray*}
for generic $x\in f\left( Y_{F}\right) $. From the exact sequence 
\begin{equation*}
0\rightarrow T_{t}\rightarrow f^{*}T_{X}\rightarrow t^{*}T_{\Bbb{P}
^{n}}\rightarrow 0
\end{equation*}
and the surjectivity of 
\begin{equation*}
T_{Y}\rightarrow \left( t\circ f\right) ^{*}T_{\Bbb{P}^{n}}
\end{equation*}
we have 
\begin{equation*}
f_{*}T_{t\circ f}\subseteq T_{t}
\end{equation*}
is a subspace of codimension $h=n-1-k.$

First suppose that 
\begin{equation*}
s_{T}=1.
\end{equation*}
Since $Y$ is $GL\left( n+1\right) $-equivariant, the normal vectors at
generic $x\in f\left( Y\right) $ can be chosen in $M_{x}^{d}$ and, since $%
s_{T}=1,$ we have that, for general $A_{1},\ldots ,A_{h}\in M_{x}^{1}$ and
for general $P\in M^{d-1}$, 
\begin{equation*}
P\cdot A_{1},\ldots ,P\cdot A_{h}\in M_{x}^{d}
\end{equation*}
generate 
\begin{equation*}
\left. N_{f}\right| _{x}=\frac{M_{x}^{d}}{T}.
\end{equation*}
Then for any $A\in M^{1}$ not vanishing at $x$, $P\cdot A$ generates $%
M^{d}/M_{x}^{d}$ so, for general $B_{1},\ldots ,B_{h+1}\in M^{1}$ and
general $P\in M^{d-1}$, 
\begin{equation}
P\cdot B_{1},\ldots ,P\cdot B_{h+1}\in M^{d}  \label{vectors}
\end{equation}
generate 
\begin{equation*}
\frac{M^{d}}{T}.
\end{equation*}

Referring to $\left( \ref{numap}\right) $ let 
\begin{equation*}
\nu ^{\prime }:M^{1}\otimes \mathcal{O}_{Y_{F}}\rightarrow N_{\iota _{F}}
\end{equation*}
be the generically surjective map defined by 
\begin{equation*}
\nu ^{\prime }\left( B\right) =\nu \left( P\cdot B\right) .
\end{equation*}
The pointwise kernels of $\nu ^{\prime }$ define a rational map 
\begin{equation*}
\mu :Y_{F}\rightarrow Gr\left( k+1,M^{1}\right) .
\end{equation*}
Under the perfect pairing 
\begin{equation*}
\bigwedge\nolimits^{h+1}M^{1}\otimes
\bigwedge\nolimits^{k+1}M^{1}\rightarrow \bigwedge\nolimits^{n+1}M^{1}=\Bbb{C%
}
\end{equation*}
induced by wedge product, we have a natural isomorphism 
\begin{equation*}
\bigwedge\nolimits^{h+1}M^{1}=\left( \bigwedge\nolimits^{k+1}M^{1}\right)
^{\vee },
\end{equation*}
By $\left( \ref{cond}\right) ,$ the kernel of the linear map 
\begin{eqnarray}
\bigwedge\nolimits^{h+1}M^{1} &\rightarrow &\mathrm{Hom}\left( \left.
f^{*}\omega _{\Bbb{P}^{n}\times S}\right| _{Y_{F}},\ \omega _{Y_{F}}\right)
\label{intmap} \\
B_{1}\wedge \ldots \wedge B_{h+1} &\mapsto &\left( \eta \mapsto \left\langle
\left. P\cdot B_{1}\wedge \ldots \wedge P\cdot B_{h+1}\right| \eta
\right\rangle \right)  \notag
\end{eqnarray}
consists in those elements of $\left( \bigwedge\nolimits^{k+1}M^{1}\right)
^{\vee }$ which vanish identically on $\mu \left( Y_{F}\right) $. But the
image of $\left( \ref{intmap}\right) $ actually lies in 
\begin{equation*}
Hom\left( \left. f^{*}\omega _{\Bbb{P}^{n}\times S}\right| _{Y_{F}},\ \omega
_{Y_{F}}\left( -\left( d-1\right) \right) \right) =\omega _{Y_{F}}\left(
n+2-d\right)
\end{equation*}
since all the vectors $P\cdot B_{1},\ldots ,P\cdot B_{h+1}$ are tangent to $%
X $ at points at which $P$ is zero. Thus 
\begin{equation}
h^{0}\left( \omega _{Y_{F}}\left( n+2-d\right) \right) \geq \left( 
\begin{array}{c}
n+1 \\ 
k+1
\end{array}
\right) -r  \label{genusbound}
\end{equation}
where 
\begin{equation}
r=\dim \left\{ h\in H^{0}\left( \mathcal{O}_{Gr\left( k+1,n+1\right) }\left(
1\right) \right) :\left. h\right| _{\mu \left( Y_{F}\right) }=0\right\} .
\label{r}
\end{equation}

We conclude the following.

\begin{lemma}
\label{s=1}Suppose that, for generic $x\in f\left( Y_{F}\right) $, 
\begin{equation*}
s_{T_{t\circ f}}=1.
\end{equation*}
Then, for $r$ as in $\left( \ref{r}\right) $, 
\begin{equation*}
h^{0}\left( \omega _{Y_{F}}\left( n+2-d\right) \right) \geq \left( 
\begin{array}{c}
n+1 \\ 
k+1
\end{array}
\right) -r.
\end{equation*}
\end{lemma}

\subsection{Versal sub-families, first reduction}

We next return to the case of arbitrary $s_{T}$. From the global generation
of $\left( \ref{gen}\right) $ we conclude that the map 
\begin{equation}
H^{0}\left( \left( t\circ \rho \right) ^{*}\bigwedge\nolimits^{h}M_{\Bbb{P}
^{n}}^{r}\right) \otimes \mathcal{O}_{\Bbb{P}^{n}}\left( a^{\prime }\right)
\rightarrow \left( \bigwedge\nolimits^{h}\frac{t^{*}M_{\Bbb{P}^{n}}^{r}}{T}
\right) \otimes \mathcal{O}_{\Bbb{P}^{n}}\left( a^{\prime }\right)
\label{gs}
\end{equation}
is (generically) surjective whenever $a^{\prime }\geq s_{T}$.

\begin{lemma}
\label{GZ}Let $Y/S^{\prime }$ be a versal sub-family of $k$-folds. Suppose
that, for generic $F\in S$, 
\begin{equation*}
\left. \left( \bigwedge\nolimits^{n-1-k}t^{*}M_{\Bbb{P}^{n}}^{d}\right)
\otimes \omega _{X}\left( a\right) \right| _{X_{F}}
\end{equation*}
is generated by global sections. Then 
\begin{equation}
H^{0}\left( \omega _{Y_{F}}\left( a\right) \right) \neq 0.  \label{cm}
\end{equation}
\end{lemma}

\begin{proof}
At generic $\left( x,F\right) \in f\left( Y\right) $, there is a vector 
\begin{equation*}
\nu \in H^{0}\left( \left. \bigwedge\nolimits^{n-1-k}t^{*}M_{\Bbb{P}
^{n}}^{d}\otimes \omega _{X}\left( a\right) \right| _{X_{F}}\right)
\end{equation*}
such that 
\begin{equation}
\left\langle \left. \nu \wedge w\wedge z\right| \left( \left. \omega
_{X}\left( a\right) \right| _{\left( x,F\right) }\right) \right\rangle \neq 0
\label{dil}
\end{equation}
where 
\begin{eqnarray*}
w &\in &\bigwedge\nolimits^{N}\left. T_{S}\right| _{F} \\
z &\in &\bigwedge\nolimits^{k}\left. f_{*}T_{Y_{F}}\right| _{\left(
x,F\right) }.
\end{eqnarray*}
But this means that 
\begin{equation*}
\left\langle \left. \nu \wedge w\right| \left( \left. \omega _{X}\left(
a\right) \right| _{Y_{F}}\right) \right\rangle
\end{equation*}
gives a non-zero element of 
\begin{equation*}
H^{0}\left( \omega _{Y_{F}}\left( a\right) \right) .
\end{equation*}
\end{proof}

We conclude:

\begin{corollary}
\label{corGZ}Let 
\begin{equation*}
f=\left( \tilde{x},\tilde{F}\right) :Y\rightarrow X
\end{equation*}
be a versal sub-family of $k$-folds for which $\left( \ref{fund}\right) $
holds.

i) 
\begin{equation*}
d+a-\left( n+1\right) <s_{T_{t\circ f}}\leq h=\left( n-1\right) -k.
\end{equation*}

ii) If $s_{T_{t\circ f}}^{\prime }+1<s_{T_{t\circ f}}$, then the map $f$
lifts to a $GL\left( n+1\right) $-equivariant map 
\begin{equation*}
g=\left( \tilde{l},\tilde{x},\tilde{F}\right) :Y\rightarrow \Delta 
\end{equation*}
such that the corank of 
\begin{equation*}
\varepsilon =e^{d}\circ g_{*}:T_{t\circ \rho \circ g}\rightarrow T_{t\circ
\rho }\rightarrow \pi ^{*}E_{\Bbb{P}^{n}}^{d}
\end{equation*}
is 
\begin{equation*}
s_{T_{t\circ f}}
\end{equation*}
and, for $T=T_{t\circ \rho \circ g}$ at a general point of $Y$ we have 
\begin{equation*}
\dim \frac{M_{l_{T}}^{d}}{T\cap M_{l_{T}}^{d}}=\dim \frac{T+M_{l_{T}}^{d}}{T}%
=rank\left( \mu _{T,s_{T}^{\prime }}\right) -s_{T}^{\prime }.
\end{equation*}
\end{corollary}

\begin{proof}
i) See $\left( \ref{gs}\right) $.

ii) Define the lifting $g$ by 
\begin{eqnarray*}
\tilde{l} &:&Y\rightarrow G \\
y &\mapsto &l_{T_{t\circ \rho \circ g}.}
\end{eqnarray*}
and use Lemma \ref{goodlem}iii).
\end{proof}

\begin{remark}
\label{list'}For 
\begin{equation*}
T=T_{t\circ f}
\end{equation*}
we recall the inequality 
\begin{equation*}
d+a-\left( n+1\right) <s_{T}\leq n-1-k
\end{equation*}
coming from Corollary \ref{corGZ}. Referring to remark \ref{list} 
\begin{equation}
\left( s_{T}-s_{T}^{\prime }\right) \geq rank\left( \mu _{T,s_{T}^{\prime
}}\right) \ \Leftrightarrow \ \left( s_{T}-s_{T}^{\prime }\right) \geq \frac{%
n-1-k}{2}  \label{equiv}
\end{equation}
and 
\begin{equation*}
s_{T}+s_{T}^{\prime }\leq \left( n-1-k\right) .
\end{equation*}
So 
\begin{eqnarray*}
s_{T}-s_{T}^{\prime } &=&2s_{T}-\left( s_{T}+s_{T}^{\prime }\right)  \\
&\geq &2s_{T}-\left( n-1-k\right)  \\
&\geq &2\left( d+a-n\right) -\left( n-1-k\right) .
\end{eqnarray*}
Combining these inequalities we have that 
\begin{equation*}
s_{T}-s_{T}^{\prime }\geq rank\left( \mu _{T,s_{T}^{\prime }}\right) 
\end{equation*}
whenever 
\begin{equation}
d+a\geq \frac{7n-3-3k}{4}.  \label{ineq1}
\end{equation}
Also 
\begin{equation*}
\left( s_{T}-s_{T}^{\prime }\right) \geq 2
\end{equation*}
whenever 
\begin{equation*}
s_{T}\geq \frac{n-1-k}{2}+1
\end{equation*}
which is insured by 
\begin{equation*}
d+a-n\geq \frac{n-1-k}{2}+1
\end{equation*}
that is, by 
\begin{equation}
d+a\geq \frac{3n+1-k}{2}.  \label{ineq 2}
\end{equation}
\end{remark}

Thus the rest of this paper is devoted to proving the assertion of the Main
Theorem in the case in which:

\begin{description}
\item[Condition 1]  The map $f$ lifts to a $GL\left( n+1\right) $
-equivariant map 
\begin{equation*}
g=\left( \tilde{l},\tilde{x},\tilde{F}\right) :Y\rightarrow \Delta .
\end{equation*}
such that 
\begin{equation*}
M_{\tilde{l}}^{d}\subseteq T_{t\circ f}+image\left( \mu _{T_{t\circ
f},s_{T_{t\circ f}}^{\prime }}\right) .
\end{equation*}

\item[Condition 2]  
\begin{equation*}
\left( s_{T_{t\circ f}}-s_{T_{t\circ f}}^{\prime }\right) \geq rank\left(
\mu _{T_{t\circ f},s_{T_{t\circ f}}^{\prime }}\right) .
\end{equation*}
(Notice that $\left( \ref{ineq1}\right) $ implies Condition 2 and $\left( 
\ref{ineq 2}\right) $ implies Condition 1.) The critical point in what
follows will be that, under these assumptions, the map 
\begin{equation*}
g=\left( \tilde{l},\tilde{x},\tilde{F}\right) :Y\rightarrow \Delta 
\end{equation*}
generically has the property 
\begin{equation*}
\left( p\circ \pi \right) M_{G}^{d}\cap f_{*}\left( \left. T_{Y}\right|
_{y}\right) \subseteq f_{*}T_{\pi \circ g}
\end{equation*}
where as above $T_{\pi \circ g}$ is the tangent space to the fibers of 
\begin{equation*}
\pi \circ g=\left( \tilde{l},\tilde{x}\right) .
\end{equation*}
Establishing this last property is perhaps the deepest part of the proof. It
is here that the Lie bracket computation introduced by Voisin is the central
ingredient.
\end{description}

\section{Voisin's bracket map}

\subsection{The contact isomorphism}

We define 
\begin{equation*}
M_{kG}^{r}\subseteq G\times M^{r}
\end{equation*}
as the sub-bundle whose fiber at $l$ is given by those forms which vanish to
order $k$ along $l$. Let 
\begin{equation*}
F\left( l^{\prime }\right) 
\end{equation*}
be any local section of the bundle $M_{G}^{d}$ near some fixed $l\in G$.
Differentiating the equation
\begin{equation*}
F\left( l^{\prime }\right) \equiv 0,
\end{equation*}
the restriction of 
\begin{equation*}
-\left. \frac{\partial F\left( l^{\prime }\right) }{\partial l^{\prime }}%
\right| _{l^{\prime }=l}
\end{equation*}
to the line $l$ itself depends only on the value of $F\left( l\right) $. So
it gives a well defined bilinear map 
\begin{equation}
M_{G}^{d}\times T_{G}\rightarrow E^{d}.  \label{map}
\end{equation}
Alternatively we can understand this pairing by viewing 
\begin{equation*}
l^{\prime }=l_{u}^{\prime }:l\rightarrow \Bbb{P}^{n}
\end{equation*}
as a family of maps from the fixed projective line $l$ with parameter $u$
(which is given by the identity map at $u=0$) and differentiating the
relation 
\begin{equation*}
F\left( l_{u}^{\prime }\right) \circ l_{u}^{\prime }=0.
\end{equation*}
We then see that 
\begin{equation}
\underset{u\rightarrow 0}{\lim }\frac{F\left( l\right) \circ l^{\prime
}\left( u\right) }{u}=-\left. \frac{\partial F\left( l_{u}^{\prime }\right) 
}{\partial u}\right| _{u=0}\circ l.  \label{lim}
\end{equation}

The map $\left( \ref{map}\right) $ factors through 
\begin{equation*}
\frac{M_{G}^{d}}{M_{2G}^{d}}\otimes T_{G}\rightarrow E^{d}
\end{equation*}
and the induced map 
\begin{equation*}
\frac{M_{G}^{d}}{M_{2G}^{d}}\rightarrow Hom\left( T_{G},E^{d}\right) 
\end{equation*}
is injective. If we restrict the lifted map to 
\begin{equation*}
\frac{M_{G}^{d}}{M_{2G}^{d}}\rightarrow Hom\left( T_{q},p^{*}E^{d}\right) 
\end{equation*}
it remains injective and the image lies the sub-bundle 
\begin{equation*}
E_{\Bbb{P}^{n}}^{d}=\left\{ \left( \alpha ,\left( l,x\right) \right) \in
p^{*}E^{d}:\alpha \left( x\right) =0\right\} .
\end{equation*}
So, by dimension, the induced map 
\begin{equation}
\frac{M_{G}^{d}}{M_{2G}^{d}}\rightarrow Hom\left( T_{q},E_{\Bbb{P}%
^{n}}^{d}\right)   \label{contact'}
\end{equation}
is an isomorphism which we call the \textit{contact isomorphism}. We rewrite 
$\left( \ref{contact'}\right) $ as a multiplication 
\begin{equation}
\bullet :\frac{M_{G}^{d}}{M_{2G}^{d}}\otimes T_{q}\rightarrow E_{\Bbb{P}%
^{n}}^{d}.  \label{mult'}
\end{equation}
Since 
\begin{equation*}
E_{\Bbb{P}^{n}}^{r}=E^{r-1}\otimes \left( \frac{U}{\Bbb{C}\cdot q}\right)
^{\vee }
\end{equation*}
and 
\begin{equation*}
T_{q}=Hom\left( \frac{U}{\Bbb{C}\cdot q},\frac{\left( M^{1}\right) ^{\vee }}{%
U}\right) ,
\end{equation*}
we can rewrite $\left( \ref{contact'}\right) $ as 
\begin{equation}
\frac{M_{G}^{d}}{M_{2G}^{d}}\rightarrow U^{\bot }\otimes E^{d-1}
\label{contact}
\end{equation}
and $\left( \ref{mult'}\right) $ as 
\begin{equation}
\bullet :\frac{M_{G}^{d}}{M_{2G}^{d}}\otimes \left( U^{\bot }\right) ^{\vee
}\rightarrow E^{d-1}  \label{mult}
\end{equation}

\subsection{The vertical contact distribution}

Referring to the exact sequences 
\begin{equation*}
0\rightarrow M_{1}^{d}=T_{t}\rightarrow T_{X}\rightarrow t^{*}T_{\Bbb{P}
^{n}}\rightarrow 0
\end{equation*}
we next wish to examine the lift of the distribution 
\begin{equation}
\left( p\circ \pi \right) ^{*}M_{G}^{d}\subset \pi ^{*}M_{\Bbb{P}
^{n}}^{d}\subset T_{\Bbb{\Delta }}  \label{dist}
\end{equation}
under maps 
\begin{equation*}
g:Y\rightarrow \Delta
\end{equation*}
such that $f=\rho \circ g$ is a versal sub-family.

First we can consider an element 
\begin{equation*}
F\in M^{d}
\end{equation*}
as 
\begin{eqnarray*}
d_{\Bbb{C}^{n+1}}F &=&\sum\nolimits_{i=0}^{n}\frac{\partial F}{\partial X_{i}%
}\otimes dX_{i} \\
&\in &H^{0}\left( \mathcal{O}_{\Bbb{P}^{n}}\left( d-1\right) \right) \otimes
H^{0}\left( \mathcal{O}_{\Bbb{P}^{n}}\left( 1\right) \right) .
\end{eqnarray*}
If we are evaluating $d_{\Bbb{C}^{n+1}}F$ at points of a line $l\subseteq
X_{F}$ then 
\begin{equation*}
\left. d_{\Bbb{C}^{n+1}}F\right| _{l}\in \left. E^{d-1}\otimes U^{\bot
}\right| _{l}
\end{equation*}
so we get a map 
\begin{equation*}
\delta _{\Bbb{C}^{n+1}}:\left( p\circ \pi \right) ^{*}M_{G}^{d}\rightarrow
E^{d-1}\otimes U^{\bot }.
\end{equation*}
The associated map 
\begin{equation*}
\pi ^{*}\left( U^{\bot }\right) ^{\vee }\otimes \left( p\circ \pi \right)
^{*}M_{G}^{d}\rightarrow E^{d-1}
\end{equation*}
is just the multiplication map given in $\left( \ref{mult}\right) $, that
is, 
\begin{equation*}
\upsilon \otimes F\mapsto \upsilon \bullet F.
\end{equation*}

Suppose now we are given a $GL\left( n+1\right) $-equivariant map 
\begin{eqnarray*}
g &:&Y\rightarrow \Delta  \\
y &\mapsto &\left( \tilde{l},\tilde{x},\tilde{F}\right) 
\end{eqnarray*}
which is an immersion at $y\mapsto \left( l,x,F\right) $. Recalling $\left( 
\ref{sum}\right) $ inside
\begin{equation*}
T_{q\circ \pi \circ g}\subseteq \left( q\circ \pi \circ g\right) ^{*}M_{\Bbb{%
\ P}^{n}}^{d}.
\end{equation*}
define the distribution
\begin{equation}
T^{\prime }=T_{q\circ \pi \circ g}\cap \left( p\circ \pi \circ g\right)
^{*}M_{G}^{d}.  \label{condition}
\end{equation}
Recalling $\left( \ref{mult'}\right) $ and the isomorphism 
\begin{equation*}
T_{q}=Hom\left( \frac{U}{\Bbb{C}\cdot q},\frac{\left( M^{1}\right) ^{\vee }}{%
U}\right) ,
\end{equation*}
we have the composition map 
\begin{eqnarray*}
T_{q\circ \pi \circ g}\overset{g_{*}}{\longrightarrow }T_{\pi }\oplus \pi %
^{*}T_{q} &\rightarrow &T_{\pi }\overset{e^{d}}{\longrightarrow }\pi ^{*}E_{%
\Bbb{P}^{n}}^{d}=\left( p\circ \pi \right) ^{*}E^{d-1}\otimes \left( \frac{U%
}{\Bbb{C}\cdot q}\right) ^{\vee } \\
\tau  &\mapsto &\tilde{F}_{*}\left( \tau \right) \oplus \tilde{l}_{*}\left(
\tau \right) \mapsto \tilde{F}_{*}\left( \tau \right) \mapsto \left. \tilde{F%
}_{*}\left( \tau \right) \right| _{l}
\end{eqnarray*}
which we denote as
\begin{equation*}
\varepsilon :T_{q\circ \pi \circ g}\rightarrow \pi ^{*}E_{\Bbb{P}^{n}}^{d}.
\end{equation*}
We are finally ready to present the second-order tool and its geometric
interpretation which are the central ingredients of everything which follows.

\begin{lemma}
\label{sff}For vector fields $\tau $,$\tau ^{\prime }$ in $T^{\prime }$,
\begin{equation*}
\varepsilon \left( \left[ \tau ,\tau ^{\prime }\right] \right) =\tilde{l}%
_{*}\left( \tau ^{\prime }\right) \bullet \tilde{F}_{*}\left( \tau \right) -%
\tilde{l}_{*}\left( \tau \right) \bullet \tilde{F}_{*}\left( \tau ^{\prime
}\right) .
\end{equation*}
\end{lemma}

\begin{proof}
Once the machinery is set up, the verification of the Lemma is
straightforward and most easily checked by an elementarty computation in
normalized local coordinates. Let $\left\{ y_{j}\right\} $ be local
coordinates for $Y$ such that 
\begin{eqnarray*}
\tilde{l}\left( 0\right)  &=&l \\
\tilde{F}\left( 0\right)  &=&F.
\end{eqnarray*}
Suppose that 
\begin{eqnarray*}
x &=&\left( 1,0,\ldots ,0\right)  \\
l &=&\left\{ X_{2}=\ldots =X_{n}=0\right\} .
\end{eqnarray*}
Then
\begin{equation*}
X_{i}=b_{i}X_{1},\ i=2,\ldots ,n,
\end{equation*}
become the local coordinates for a small neighborhood of $l$ in $G_{x}$.
Letting $J$ denote a multi-index for the variables $\left\{ y_{j}\right\} $
we write 
\begin{equation*}
\tilde{l}\left( \left\{ y_{j}\right\} \right) \mapsto \left\{ b_{i}\left(
\left\{ y_{j}\right\} \right) \right\} =\left\{ \sum\nolimits_{J}b_{i,J}%
\cdot y^{J}\right\} .
\end{equation*}
Also if $I$ denotes multi-indices for the variables $\left\{ X_{2},\ldots
,X_{n}\right\} $ we write 
\begin{equation*}
\tilde{F}\left( \left\{ y_{j}\right\} \right) =\sum\nolimits_{I,J}y^{J}%
\tilde{F}_{I;J}X^{I}
\end{equation*}
where 
\begin{equation*}
\tilde{F}_{I;J}=\tilde{F}_{I;J}\left( X_{0},X_{1}\right) 
\end{equation*}
is a homogeneous form of degree $d-\left| I\right| $. The condition that
\begin{equation*}
\sum\nolimits_{k}a_{k}^{\prime }\frac{\partial }{\partial y_{k}}=\tau
^{\prime }\in T^{\prime }
\end{equation*}
is the equation, for all $\left\{ y_{j}\right\} $ that
\begin{equation}
\sum\nolimits_{k,I}a_{k}^{\prime }\frac{\partial \sum\nolimits_{I}y^{J}%
\tilde{F}_{I;J}}{\partial y_{k}}b^{I}\left( \left\{ y_{j}\right\} \right)
X_{1}^{\left| I\right| }=0.  \label{expression}
\end{equation}
Applying $\frac{\partial }{\partial y_{j}}$ to $\left( \ref{expression}%
\right) $ and then setting $\left\{ y_{j}\right\} =0$ we obtain 
\begin{equation*}
\frac{\partial a_{k}^{\prime }}{\partial y_{j}}\tilde{F}_{0;k}+\sum%
\nolimits_{k}a_{k}^{\prime }\left( 1+\delta _{jk}\right) \tilde{F}%
_{0;jk}+\sum\nolimits_{i=2}^{n}\sum\nolimits_{k}a_{k}^{\prime }\frac{%
\partial b_{i}}{\partial y_{j}}\left( 0\right) \cdot \tilde{F}_{i;k}=0.
\end{equation*}
So for
\begin{equation*}
\sum\nolimits_{j}a_{j}\frac{\partial }{\partial y_{j}}=\tau \in T^{\prime },
\end{equation*}
we can compute
\begin{equation*}
\left[ \tau ,\tau ^{\prime }\right] \tilde{F}=\sum\nolimits_{j}a_{j}\frac{%
\partial a_{k}^{\prime }}{\partial y_{j}}\tilde{F}_{0;k}-\sum%
\nolimits_{j}a_{j}^{\prime }\frac{\partial a_{k}}{\partial y_{j}}\tilde{F}%
_{0;k}
\end{equation*}
at $\left\{ y_{j}\right\} =0$ as
\begin{equation}
-\sum\nolimits_{i=2}^{n}\sum\nolimits_{j,k}a_{j}\frac{\partial b_{i}}{%
\partial y_{j}}\left( 0\right) \cdot a_{k}^{\prime }\tilde{F}%
_{i;k}+\sum\nolimits_{i=2}^{n}\sum\nolimits_{j,k}a_{j}^{\prime }\frac{%
\partial b_{i}}{\partial y_{j}}\left( 0\right) \cdot a_{k}\tilde{F}_{i;k}.
\label{EXP}
\end{equation}
Rewriting $\left( \ref{EXP}\right) $ as 
\begin{equation}
\sum\nolimits_{j,k}\sum\nolimits_{i=2}^{n}\left( a_{k}^{\prime }\frac{%
\partial b_{i}}{\partial y_{k}}\left( 0\right) \cdot a_{j}\frac{\partial ^{2}%
\tilde{F}}{\partial X_{i}\partial y_{j}}\left( 0\right) -a_{j}\frac{\partial %
b_{i}}{\partial y_{j}}\left( 0\right) \cdot a_{k}^{\prime }\frac{\partial %
^{2}\tilde{F}}{\partial X_{i}\partial y_{k}}\left( 0\right) \right) 
\label{obst'}
\end{equation}
and restricting to $\tilde{l}\left( 0\right) $ we see that we obtain exactly
\begin{equation*}
\tilde{l}_{*}\left( \tau ^{\prime }\right) \bullet \tilde{F}_{*}\left( \tau
\right) -\tilde{l}_{*}\left( \tau \right) \bullet \tilde{F}_{*}\left( \tau
^{\prime }\right) .
\end{equation*}
\end{proof}

\begin{corollary}
\label{corsff}Suppose that Conditions 1 and 2 above hold. Let $T^{\prime }$
be as in the lemma. Then at generic $y=\left( l,x,F\right) \in Y$, 
\begin{equation*}
\tilde{l}_{*}\left( \left. T^{\prime }\right| _{y}\cap M_{2l}^{d}\right) =0.
\end{equation*}
\end{corollary}

\begin{proof}
Recall that 
\begin{equation*}
codim\left( Y,X\right) =\dim \frac{M_{x}^{d}}{\left. T_{q\circ \pi \circ
g}\right| _{y}}=n-1-k.
\end{equation*}
Also by Corollary \ref{corGZ}ii) we have for $T=\left. T_{q\circ \pi \circ
g}\right| _{y}$ that 
\begin{equation*}
\dim \left( \varepsilon \left( \left. T_{q\circ \pi \circ g}\right|
_{y}\right) \right) =d-s_{T}.
\end{equation*}
Suppose 
\begin{equation*}
\tau ^{\prime }\in \left. T_{q\circ \pi \circ g}\right| _{y}\cap \tilde{l}
^{*}M_{2l}^{d}.
\end{equation*}
Then for every $\tau \in T_{y}^{\prime }$ the element 
\begin{equation*}
\tilde{l}_{*}\left( \tau ^{\prime }\right) \bullet \tilde{F}_{*}\left( \tau
\right) -\tilde{l}_{*}\left( \tau \right) \bullet \tilde{F}_{*}\left( \tau
^{\prime }\right) =\tilde{l}_{*}\left( \tau ^{\prime }\right) \bullet \tilde{
F}_{*}\left( \tau \right)
\end{equation*}
lies in the image of $\varepsilon $. If $\tilde{l}_{*}\left( \tau ^{\prime
}\right) \neq 0$, then the map 
\begin{equation*}
\tilde{\psi}:\frac{M_{l}^{d}}{M_{2l}^{d}}=\left. E_{\Bbb{P}^{n}}^{d}\otimes
T_{q}^{\vee }\right| _{\left( l,x\right) }\overset{\tilde{l}_{*}\left( \tau
^{\prime }\right) }{\longrightarrow }\left. E_{\Bbb{P}^{n}}^{d}\right|
_{\left( l,x\right) }
\end{equation*}
is surjective and so 
\begin{equation*}
\dim \left( \ker \tilde{\psi}\right) =\left( n-1\right) d-d.
\end{equation*}
Let 
\begin{equation*}
\psi =\left. \tilde{\psi}\right| _{T_{y}^{\prime }/T_{y}^{\prime }\cap
M_{2l}^{d}}.
\end{equation*}
But 
\begin{equation*}
d-s_{T}\geq \dim \psi \left( \frac{T_{y}^{\prime }}{T_{y}^{\prime }\cap
M_{2l}^{d}}\right) .
\end{equation*}
So 
\begin{eqnarray*}
d-s_{T} &\geq &\dim \left( \frac{T_{y}^{\prime }}{T_{y}^{\prime }\cap
M_{2l}^{d}}\right) -\left( \left( n-1\right) d-d\right) \\
\left( n-1\right) d-s_{T} &\geq &\dim \left( \frac{T_{y}^{\prime }}{
T_{y}^{\prime }\cap M_{2l}^{d}}\right) .
\end{eqnarray*}
On the other hand from Corollary \ref{corGZ}ii) we have for $T=\left.
T_{q\circ \pi \circ g}\right| _{y}$ that 
\begin{equation}
\dim \frac{M_{l}^{d}}{T_{y}^{\prime }}=rank\left( \mu _{T,s_{T}^{\prime
}}\right) -s_{T}^{\prime }.  \label{previneq}
\end{equation}
So 
\begin{equation*}
\dim \frac{T_{y}^{\prime }}{T_{y}^{\prime }\cap M_{2l}^{d}}=\dim \frac{
T_{y}^{\prime }+M_{2l}^{d}}{M_{2l}^{d}}=\left( n-1\right) d-\dim \frac{
M_{l}^{d}}{T_{y}^{\prime }+M_{2l}^{d}}
\end{equation*}
and therefore by $\left( \ref{previneq}\right) $ we have 
\begin{equation*}
\dim \frac{T_{y}^{\prime }}{T_{y}^{\prime }\cap M_{2l}^{d}}\geq \left(
n-1\right) d-\left( rank\left( \mu _{T,s_{T}^{\prime }}\right)
-s_{T}^{\prime }\right) .
\end{equation*}
Thus 
\begin{equation*}
\left( n-1\right) d-s_{T}\geq \left( n-1\right) d-rank\left( \mu
_{T,s_{T}^{\prime }}\right) +s_{T}^{\prime }
\end{equation*}
that is 
\begin{equation*}
rank\left( \mu _{T,s_{T}^{\prime }}\right) \geq s_{T}+s_{T}^{\prime }.
\end{equation*}
But by $\left( \ref{ineq}\right) $%
\begin{equation*}
\left( n-1-k\right) -\left( s_{T}-s_{T}^{\prime }\right) =\dim \left(
image\left( \mu _{T,s_{T}^{\prime }}\right) \right)
\end{equation*}
and so we have 
\begin{equation}
n-1-k\geq 2s_{T}  \label{weird}
\end{equation}
which contradicts $\left( \ref{equiv}\right) $ unless 
\begin{equation*}
s_{T}^{\prime }=0
\end{equation*}
But then 
\begin{equation*}
s_{T}=n-1-k
\end{equation*}
which again contradicts $\left( \ref{weird}\right) $.
\end{proof}

\subsection{A critical lemma from linear algebra}

We will need the following linear algebra computation. By Corollary \ref
{corsff} the map 
\begin{equation}
\tilde{l}_{*}:T^{\prime }\rightarrow \left. T_{q}\right| _{\left( l,x\right)
}  \label{lstar}
\end{equation}
induces 
\begin{equation*}
\left. \frac{T^{\prime }}{T^{\prime }\cap M_{2G}^{d}}\right| _{y}\rightarrow
\left. T_{q}\right| _{\left( l,x\right) }
\end{equation*}
where, for notational simplicity, we denote $\left( p\circ \pi \right)
^{*}M_{2G}^{d}$ simply as $M_{2G}^{d}$. Since 
\begin{equation*}
\left. \frac{T^{\prime }}{T^{\prime }\cap M_{2G}^{d}}\right| _{y}\cong
\left. \frac{T^{\prime }+M_{2G}^{d}}{M_{2G}^{d}}\right| _{y}\subseteq \left. 
\frac{M_{G}^{d}}{M_{2G}^{d}}\right| _{y}
\end{equation*}
and, by $\left( \ref{contact'}\right) $, 
\begin{equation*}
\left. \frac{M_{G}^{d}}{M_{2G}^{d}}\right| _{y}\cong \left. \pi ^{*}\left(
E_{\Bbb{P}^{n}}^{d}\otimes T_{q}^{\vee }\right) \right| _{\left(
l,x,F\right) },
\end{equation*}
we can extend $\tilde{l}_{*}$ to a map of the same rank 
\begin{equation*}
\varphi :\left. E_{\Bbb{P}^{n}}^{d}\otimes T_{q}^{\vee }\right| _{\left(
l,x\right) }\rightarrow \left. T_{q}\right| _{\left( l,x\right) }.
\end{equation*}
Now let 
\begin{eqnarray*}
K^{\vee } &=&\left. T_{q}\right| _{x} \\
W &=&\left. E_{\Bbb{P}^{n}}^{d}\right| _{\left( l,x\right) }.
\end{eqnarray*}

On the other hand, considering 
\begin{equation*}
T^{\prime }\subseteq \left. \left( p\circ \pi \right) ^{*}M_{G}^{d}\right|
_{y}
\end{equation*}
we have 
\begin{equation*}
\vartheta :\left. T^{\prime }\right| _{y}\rightarrow \left. \frac{M_{G}^{d}}{
M_{2G}^{d}}\right| _{l}=W\otimes K
\end{equation*}
whose image we denote by $H$. The standard map 
\begin{equation}
\begin{array}{c}
\psi :\bigwedge\nolimits^{2}\left( W\otimes K\right) \rightarrow W \\ 
A\wedge B\mapsto \left\langle \left. \vartheta \left( A\right) \right|
\varphi \left( B\right) \right\rangle -\left\langle \left. \vartheta \left(
B\right) \right| \varphi \left( A\right) \right\rangle
\end{array}
\label{wedge''}
\end{equation}
restricts to a map 
\begin{equation*}
\bigwedge\nolimits^{2}H\rightarrow W
\end{equation*}
under which 
\begin{equation*}
\vartheta \left( \tau \right) \wedge \vartheta \left( \tau ^{\prime }\right)
\mapsto \varepsilon \left( \left[ \tau ,\tau ^{\prime }\right] \right) .
\end{equation*}

The needed linear algebra result is then:

\begin{lemma}
\label{la}i) Let 
\begin{equation*}
H\subseteq W\otimes K
\end{equation*}
be a subspace of codimension $c$ and 
\begin{equation*}
\varphi :H\rightarrow K^{\vee }
\end{equation*}
be any linear map and let 
\begin{equation*}
\psi :H\wedge H\rightarrow W
\end{equation*}
be defined as in $\left( \ref{wedge''}\right) .$ Let $J=\psi \left( H\wedge
H\right) $. Then 
\begin{equation*}
\dim \frac{W}{J}\leq c+1
\end{equation*}
and, if equality holds, $\varphi $ factors as a composition 
\begin{equation*}
H\rightarrow \frac{W\otimes K}{J^{\prime }\otimes K}\rightarrow K^{\vee }
\end{equation*}
where $J^{\prime }\supseteq J$ is a hyperplane in $W$, and the image of the
associated morphism 
\begin{equation*}
\frac{W}{J^{\prime }}\rightarrow K^{\vee }\otimes K^{\vee }
\end{equation*}
lies in 
\begin{equation*}
Sym^{2}\left( K^{\vee }\right) .
\end{equation*}
ii) If $c>0$ in i) then 
\begin{equation*}
\dim \frac{W}{J}\leq c.
\end{equation*}
\end{lemma}

\begin{proof}
i) Pick a complementary subspace $J^{\bot }$ to $J$ in $W$ and a basis $%
\left\{ w_{j}\right\} $ for $W$ which is compatible with the decomposition 
\begin{equation*}
W=J\oplus J^{\bot }.
\end{equation*}
Let $\left\{ k_{i}\right\} $ be a basis of $K$ and use the induced
isomorphism 
\begin{eqnarray*}
K &\rightarrow &K^{\vee } \\
k_{i} &\rightarrow &\left( k_{i^{\prime }}\mapsto \delta _{ii^{\prime
}}\right)
\end{eqnarray*}
to identify $K$ and $K^{\vee }$. Pick a complementary space $H^{\bot }$ to $%
H $ in $W\otimes K$ which has a basis consisting of monomials 
\begin{equation}
w_{j\left( h\right) }\otimes k_{i\left( h\right) }  \label{basis}
\end{equation}
for $h=1,\ldots ,c$, and extend $\varphi $ to $W$ by setting 
\begin{equation*}
\varphi \left( w_{j\left( h\right) }\otimes k_{i\left( h\right) }\right) =0.
\end{equation*}
Then for 
\begin{equation*}
\psi :\bigwedge\nolimits^{2}\left( W\otimes K\right) \rightarrow W
\end{equation*}
as in $\left( \ref{wedge''}\right) $ we have 
\begin{equation*}
\dim \left( \psi \left( H^{\bot }\otimes \left( W\otimes K\right) \right)
\right) \leq c,
\end{equation*}
so we are reduced to proving the assertion of the lemma in the case $c=0$.

In that case $\varphi $ is given by a matrix 
\begin{equation*}
\left( a_{i^{\prime }}^{j,i}\right)
\end{equation*}
and 
\begin{equation}
\psi \left( \left( w_{j}\otimes k_{i}\right) \wedge \left( w_{j^{\prime
}}\otimes k_{i^{\prime }}\right) \right) =a_{i}^{j^{\prime },i^{\prime
}}w_{j}-a_{i^{\prime }}^{j,i}w_{j^{\prime }}.  \label{formula}
\end{equation}
If $J^{\bot }\neq 0$, pick $w_{j}\in J^{\bot }$ and conclude that $%
a_{i}^{j^{\prime },i^{\prime }}=0$ unless $j=j^{\prime }$ whenever $%
w_{j^{\prime }}\in J^{\bot }$ and 
\begin{equation*}
a_{i}^{j,i^{\prime }}=a_{i^{\prime }}^{j^{\prime },i}.
\end{equation*}

ii) Assume 
\begin{equation*}
2\leq \dim \frac{W}{J}=c+1.
\end{equation*}
Then the hyperplane $J^{\prime }$ will vary non-trivially as the choice of $%
\left( \ref{basis}\right) $ varies over all possible bases.
\end{proof}

\section{Line contact for lifted $Y/S^{\prime }$}

\subsection{$Y/S^{\prime }$ lies in locus of osculating lines}

\begin{theorem}
\label{osc}Let $Y/S^{\prime }$ be a versal family of $k$-folds in $X/S$ such
that Condition 1 and Condition 2 above hold. Then, at generic $y\in Y$:

i) Either 
\begin{equation*}
s_{T_{t\circ f}}\leq 1
\end{equation*}
(in which case Lemma \ref{s=1} applies) or, referring to $\left( \ref{lstar}%
\right) $, we have
\begin{equation*}
\tilde{l}_{*}=0.
\end{equation*}
ii) If $\tilde{l}_{*}=0$ in $\left( \ref{lstar}\right) ,$ then the
distribution 
\begin{equation*}
T_{q\circ \pi \circ g}\cap M_{G}^{d}
\end{equation*}
is integrable and 
\begin{equation*}
\tilde{l}\left( y\right) 
\end{equation*}
has contact with $X_{\tilde{F}\left( y\right) }$ of order at least $d$ at $%
\tilde{x}\left( y\right) $.
\end{theorem}

\begin{proof}
i) Fix generic $y\in Y$ and let $\left( l,x,F\right) =g\left( y\right) .$
Let $\tilde{J}$ denote the image of the map $\varepsilon $ given by 
\begin{equation*}
T_{q\circ \pi \circ g}\overset{g_{*}}{\longrightarrow }T_{\pi }\oplus \pi
^{*}T_{q}\rightarrow T_{\pi }\overset{e^{d}}{\longrightarrow }\pi ^{*}E_{%
\Bbb{P}^{n}}^{d}
\end{equation*}
at $y$ and let 
\begin{equation*}
W=\left. \pi ^{*}E_{\Bbb{P}^{n}}^{d}\right| _{\left( l,x,F\right) }.
\end{equation*}
Then by Lemma \ref{goodlem}iii) 
\begin{equation*}
\dim \frac{W}{\tilde{J}}=s_{T_{t\circ f}}.
\end{equation*}
On the other hand, let 
\begin{equation*}
H
\end{equation*}
denote the image of $T^{\prime }$ under the map 
\begin{equation*}
T^{\prime }\twoheadrightarrow \left. \frac{T^{\prime }+\tilde{l}
^{*}M_{2G}^{d}}{\tilde{l}^{*}M_{2G}^{d}}\right| _{y}\subseteq \left. \frac{%
\tilde{l}^{*}M_{G}^{d}}{\tilde{l}^{*}M_{2G}^{d}}\right| _{y}=\left. \pi
^{*}\left( E_{\Bbb{P}^{n}}^{d}\otimes T_{q}\right) \right| _{\left(
l,x,F\right) }.
\end{equation*}
Then Lemma \ref{sff} and Lemma \ref{la} then imply that either $\tilde{l}
_{*}=0$ or 
\begin{equation*}
\dim \frac{W}{J}\leq 1.
\end{equation*}
But 
\begin{equation*}
\dim \frac{W}{J}\geq \dim \frac{W}{\tilde{J}}=s_{T_{t\circ f}}.
\end{equation*}
So either $s_{T_{t\circ f}}\leq 1$ or $\varphi $ must be the zero map.

ii) If $\varphi =0$ then $\psi =0$ on $T^{\prime }$ as well so that the
distribution $T^{\prime }$ is integrable. Since $\tilde{l}_{*}=0$, 
\begin{equation*}
\tilde{l}:Y\rightarrow G
\end{equation*}
is constant along $T^{\prime }$, so that 
\begin{equation*}
T^{\prime }\subseteq T_{\pi \circ g}.
\end{equation*}

Let $y$ be a general point in $Y$ with $g\left( y\right) =\left(
l,x,F\right) $. Now $\tilde{l}$, $\tilde{x}$ and 
\begin{equation*}
\left. \tilde{F}\right| _{l}
\end{equation*}
are all constant on the leaf $Y_{\left( l,x,F\right) }^{\prime }$ through $y$
which integrates $T^{\prime }$. On the other hand 
\begin{equation*}
Y_{\left( l,x\right) }=\left( \pi \circ g\right) ^{-1}\left( l,x\right)
\end{equation*}
has tangent space $\left. T_{\pi \circ g}\right| _{Y_{\left( l,x\right) }}$.
Thus 
\begin{equation*}
\dim Y_{\left( l,x\right) }-\dim Y_{\left( l,x,F\right) }^{\prime }=\mathrm{%
\ rank}\frac{T_{t\circ f}}{T^{\prime }}-\left( n-1\right) .
\end{equation*}
By Lemma \ref{goodlem} 
\begin{equation*}
rank\frac{T_{t\circ f}}{T^{\prime }}=d-s_{T_{t\circ f}}
\end{equation*}
so that 
\begin{equation*}
\dim Y_{\left( l,x\right) }-\dim Y_{\left( l,x\right) }^{\prime
}=d-s_{T_{t\circ f}}-\left( n-1\right) .
\end{equation*}
On the other hand, by Corollary \ref{corGZ} 
\begin{equation*}
d-s_{T_{t\circ f}}\leq n-a
\end{equation*}
for $a\geq 0$. Thus $a\leq 1$ and 
\begin{equation*}
\dim Y_{\left( l,x\right) }-\dim Y_{\left( l,x\right) }^{\prime }\leq 1-a.
\end{equation*}

As above we assume that $l$ is given by 
\begin{equation*}
X_{j}=0,\ j\geq 2
\end{equation*}
and 
\begin{equation*}
x=\left[ 1,0,\ldots ,0\right] .
\end{equation*}
We consider the map 
\begin{eqnarray*}
Y_{\left( l,x\right) } &\rightarrow &\left. E_{\Bbb{P}^{n}}^{d}\right|
_{\left( l,x\right) }. \\
y^{\prime } &\rightarrow &\left. \tilde{F}\left( y^{\prime }\right) \right|
_{l}
\end{eqnarray*}
The fiber of this map containing the (generically chosen) basepoint $y$ is
exactly $Y_{\left( l,x\right) }^{\prime }$. So the image is of dimension at
most 
\begin{equation*}
1-a.
\end{equation*}
On the other hand, the image is invariant under the action of the stabilizer
of $\left( l,x\right) $ in $GL\left( n+1\right) $. Letting $%
P_{X,_{0},X_{1}}^{d-1}$ denote the space of homogeneous forms of degree $d-1$
in $X_{0}$ and $X_{1}$, it is clear that the only such subsets of 
\begin{equation*}
X_{1}\cdot P_{X,_{0},X_{1}}^{d-1}
\end{equation*}
of dimension $\leq 1$ invariant under the group 
\begin{equation*}
\left\{ \left[ 
\begin{array}{cc}
\ast & * \\ 
0 & *
\end{array}
\right] \right\}
\end{equation*}
are $\left\{ 0\right\} $ and $\left\{ \Bbb{C}\cdot X_{1}^{d}\right\} $. Thus
either 
\begin{equation*}
\left. \tilde{F}\left( y\right) \right| _{l}=cX_{1}^{d}
\end{equation*}
or 
\begin{equation*}
\left. \tilde{F}\left( y\right) \right| _{l}=0.
\end{equation*}
\end{proof}

\subsection{Line osculation hierarchy}

To complete the proof of Theorem \ref{MT}, we lastly study the geometry of
the hierarchy of varieties 
\begin{equation*}
\Delta _{r}:=\left\{ \left( l,x,F\right) \in L\times _{\Bbb{P}^{n}}X:l\cdot
X_{F}\geq r\cdot x\right\} .
\end{equation*}
(We only need the cases $r=d$ and $r=d+1$ but the fundamental calculations
are the same for all $r$ so we make them in general.) We have the following
(commutative) diagram of maps: 
\begin{equation}
\begin{array}{ccccccccc}
&  &  &  & \Delta _{d+1} &  &  &  &  \\ 
&  &  &  & \cap &  &  &  &  \\ 
&  &  &  & \ldots &  &  &  &  \\ 
&  &  &  & \cap &  &  &  &  \\ 
&  &  &  & \Delta _{1}=\Delta &  &  &  &  \\ 
&  &  & \swarrow _{\pi } &  & \searrow ^{\rho } &  &  &  \\ 
&  & L &  &  &  & X &  &  \\ 
& \swarrow _{p} &  & \searrow ^{q} &  & \swarrow _{t} &  & \searrow ^{s} & 
\\ 
G &  &  &  & \Bbb{P}^{n} &  &  &  & S
\end{array}
.  \label{BD}
\end{equation}
Write 
\begin{eqnarray*}
\pi _{r} &=&\left. \pi \right| _{\Delta _{r}}:\Delta _{r}\rightarrow L \\
\rho _{r} &=&\left. \rho \right| _{\Delta _{r}}:\Delta _{r}\rightarrow X.
\end{eqnarray*}
Since the fibers of 
\begin{equation*}
\pi _{r}:\Delta _{r}\rightarrow L
\end{equation*}
are punctured vector spaces of dimension 
\begin{equation*}
N-r,
\end{equation*}
$\Delta _{r}$ is smooth and irreducible of dimension 
\begin{equation*}
N-r+2\left( n-1\right) +1
\end{equation*}
for each $r$.

\begin{lemma}
i) The map 
\begin{equation*}
\rho _{r}:\Delta _{r}\rightarrow X
\end{equation*}
is surjective as long as 
\begin{equation*}
r\leq n
\end{equation*}
and generically injective when 
\begin{equation*}
r>n.
\end{equation*}

ii) The map 
\begin{equation*}
s\circ \rho _{r}:\Delta _{r}\rightarrow S
\end{equation*}
is surjective if 
\begin{equation*}
r\leq 2\left( n-1\right) 
\end{equation*}
so that 
\begin{equation*}
\Delta _{r,F}:=\left( s\circ \rho _{r}\right) ^{-1}\left( F\right) 
\end{equation*}
is smooth for generic $F\in S^{d}$ in that case.
\end{lemma}

\begin{proof}
i) Assume 
\begin{equation*}
x=\left[ 1,0,\ldots ,0\right] \in Proj\left( \Bbb{C}\left[ X_{0},\ldots
,X_{n}\right] \right) .
\end{equation*}
The contact conditions for a line through $x$ with respect to 
\begin{equation*}
F=\sum\nolimits_{j=1}^{n}X_{0}^{d-j}F_{j}\left( X_{1},\ldots ,X_{n}\right)
\end{equation*}
become 
\begin{equation*}
\left\{ F_{1}=\ldots =F_{r-1}=0\right\} \subseteq Proj\left( \Bbb{C}\left[
X_{1},\ldots ,X_{n}\right] \right)
\end{equation*}
ii) A constant count shows that all hypersurfaces $X_{F}$ in $\Bbb{P}^{n}$
admit lines with a point of contact of order $2n-2$.
\end{proof}

Let 
\begin{equation*}
M_{\Bbb{P}^{n}}^{r}=\ker \left( M^{r}\otimes \mathcal{O}_{\Bbb{P}^{n}}%
\overset{eval.}{\longrightarrow }\mathcal{O}_{\Bbb{P}^{n}}\left( r\right)
\right)
\end{equation*}
and 
\begin{equation*}
M_{G}^{r}=\ker \left( M^{r}\otimes \mathcal{O}_{G}\longrightarrow
E^{r}\right) .
\end{equation*}
We have as spaces that 
\begin{eqnarray*}
\Delta _{d+1} &=&p^{*}M_{G}^{d}-\left\{ 0\right\} \\
\Delta &=&q^{*}M_{\Bbb{P}^{n}}^{d}-\left\{ 0\right\} .
\end{eqnarray*}

\subsection{Canonical bundles of spaces of line-osculators}

Now 
\begin{equation*}
U^{\vee }=p_{*}q^{*}\mathcal{O}_{\Bbb{P}^{n}}\left( 1\right)
\end{equation*}
so that 
\begin{equation*}
c_{1}\left( U^{\vee }\right) =\mathcal{O}_{G}\left( 1\right) .
\end{equation*}
Also 
\begin{equation*}
T_{G}=Hom\left( U,\frac{G\times \Bbb{C}^{n+1}}{U}\right)
\end{equation*}
so that 
\begin{equation*}
c_{1}\left( \omega _{G}\right) =\left( n+1\right) \cdot c_{1}\left( U\right)
.
\end{equation*}
For $E^{d}$ as defined above, we have 
\begin{equation*}
E^{d}=p_{*}q^{*}\mathcal{O}_{\Bbb{P}^{n}}\left( d\right) =Sym^{d}U^{\vee }
\end{equation*}
so that 
\begin{equation*}
c_{1}\left( E^{d}\right) =\mathcal{O}_{G}\left( \frac{d\left( d+1\right) }{2}%
\right) .
\end{equation*}
$p^{*}E^{r}$ has a distinguished line sub-bundle 
\begin{equation*}
\mathcal{L}_{L}^{r}
\end{equation*}
whose fiber at $\left( x,l\right) \in L$ is 
\begin{equation*}
H^{0}\left( \mathcal{O}_{l}\left( r\right) \left( -r\cdot x\right) \right) .
\end{equation*}
We have 
\begin{eqnarray*}
c_{1}\left( \mathcal{L}_{L}^{r}\right) &=&r\cdot c_{1}\left( \mathcal{L}%
_{L}^{1}\right) \\
&=&p^{*}\mathcal{O}_{G}\left( r\right) -q^{*}\mathcal{O}_{\Bbb{P}^{n}}\left(
r\right)
\end{eqnarray*}
For $r\leq d$, define 
\begin{equation*}
\frak{M}_{r}^{d}=\mathcal{L}_{L}^{r}\otimes p^{*}E^{d-r}\subseteq p^{*}E^{d}.
\end{equation*}
so that 
\begin{eqnarray}
c_{1}\left( \frak{M}_{r}^{d}\right) &=&p^{*}c_{1}\left( E^{d-r}\right)
+\left( d-r+1\right) \cdot c_{1}\left( \mathcal{L}_{L}^{r}\right)
\label{form} \\
&=&p^{*}\mathcal{O}_{G}\left( \frac{\left( d-r\right) \left( d-r+1\right) }{2%
}\right)  \notag \\
&&+\left( d-r+1\right) \cdot r\cdot \left( p^{*}\mathcal{O}_{G}\left(
1\right) -q^{*}\mathcal{O}_{\Bbb{P}^{n}}\left( 1\right) \right) \\
&=&p^{*}\mathcal{O}_{G}\left( \frac{\left( d+r\right) \left( d-r+1\right) }{2%
}\right) -q^{*}\mathcal{O}_{\Bbb{P}^{n}}\left( r\cdot \left( d-r+1\right)
\right) .  \notag
\end{eqnarray}

Define 
\begin{equation*}
\frak{F}_{r}^{d}=\frac{p^{*}E^{d}}{\frak{M}_{r}^{d}}
\end{equation*}
and 
\begin{equation*}
M_{r}^{d}=\ker \left( M^{d}\otimes \mathcal{O}_{L}\rightarrow \frak{F}
_{r}^{d}\right) .
\end{equation*}
Notice that 
\begin{eqnarray*}
M_{1}^{d} &=&M_{\Bbb{P}^{n}}^{d} \\
M_{d+1}^{d} &=&p^{*}M_{G}^{d}.
\end{eqnarray*}
Also 
\begin{equation*}
\Delta _{r}=M_{r}^{d}-\left\{ 0\right\}
\end{equation*}
so that in the exact tangent bundle sequence 
\begin{equation*}
0\rightarrow T_{\pi }\rightarrow T_{\Delta _{r}}\rightarrow \pi
^{*}T_{L}\rightarrow 0.
\end{equation*}
we have 
\begin{equation*}
T_{\pi }=\pi ^{*}M_{r}^{d}.
\end{equation*}
Furthermore 
\begin{eqnarray*}
c_{1}\left( M_{r}^{d}\right) &=&-c_{1}\left( \frak{F}_{r}^{d}\right) \\
&=&c_{1}\left( \frak{M}_{r}^{d}\right) -c_{1}\left( p^{*}E^{d}\right) \\
&=&-p^{*}\mathcal{O}_{G}\left( \frac{r\left( r-1\right) }{2}\right) -q^{*}%
\mathcal{O}_{\Bbb{P}^{n}}\left( r\cdot \left( d-r+1\right) \right) .
\end{eqnarray*}
Notice that this formula also holds for $r=d+1$.

Also we have 
\begin{equation*}
\det \frak{F}_{r}^{d}=p^{*}\mathcal{O}_{G}\left( \frac{r\left( r-1\right) }{2%
}\right) +q^{*}\mathcal{O}_{\Bbb{P}^{n}}\left( r\cdot \left( d-r+1\right)
\right) .
\end{equation*}
The image of the tautological section 
\begin{eqnarray*}
\Delta &\rightarrow &\left( p\circ \pi \right) ^{*}E^{d} \\
F &\mapsto &\left. F\right| _{l}
\end{eqnarray*}
in the quotient bundle 
\begin{equation*}
\frak{F}_{r}^{d}=\frac{p^{*}E^{d}}{\frak{M}_{r}^{d}}
\end{equation*}
has zero-scheme $\Delta _{r}$. and so, by adjunction, 
\begin{equation*}
\omega _{\Delta _{r}}=\left( \pi ^{*}\omega _{L/G}\right) \otimes \left( \pi
^{*}\omega _{G}\right) \otimes \det \frak{F}_{r}^{d}.
\end{equation*}

On the other hand, from the exact sequence 
\begin{equation*}
0\rightarrow \mathcal{O}_{L}\rightarrow p^{*}U\otimes q^{*}\mathcal{O}_{\Bbb{%
\ P}^{n}}\left( 1\right) \rightarrow T_{L/G}\rightarrow 0
\end{equation*}
we have 
\begin{eqnarray*}
\omega _{L/G} &=&p^{*}\det U^{\vee }\otimes q^{*}\mathcal{O}_{\Bbb{P}%
^{n}}\left( -2\right) \\
&=&p^{*}\mathcal{O}_{G}\left( 1\right) \otimes q^{*}\mathcal{O}_{\Bbb{P}
^{n}}\left( -2\right) .
\end{eqnarray*}
So 
\begin{equation}
\omega _{\Delta _{r}}=\pi ^{*}\left( p^{*}\left( \mathcal{O}_{G}\left( \frac{
r\left( r-1\right) }{2}-n\right) \right) \otimes q^{*}\mathcal{O}_{\Bbb{P}
^{n}}\left( r\left( d-r+1\right) -2\right) \right) .  \label{F1}
\end{equation}
Thus, for example, whenever 
\begin{eqnarray}
\frac{d\left( d-1\right) }{2} &\geq &n  \label{weak} \\
a+d &\geq &2  \notag
\end{eqnarray}
we conclude that 
\begin{equation*}
\omega _{\Delta _{d}}\left( a\right) =\omega _{\Delta _{d}}\otimes \left(
q\circ \pi \right) ^{*}\mathcal{O}_{\Bbb{P}^{n}}\left( a\right)
\end{equation*}
is the the pull-back of a globally generated bundle on $L$ and therefore is
globally generated. Thus by adjunction, for generic $F\in S^{d}$%
\begin{equation*}
\omega _{\Delta _{r,F}}\left( a\right)
\end{equation*}
is globally generated.

Notice that the analogous computation for $\Delta _{d+1}$ gives 
\begin{equation*}
\omega _{\Delta _{d+1}}=\left( \frac{d\left( d+1\right) }{2}-n\right) \pi
^{*}p^{*}c_{1}\left( U^{\vee }\right) +\pi ^{*}c_{1}\left( \mathcal{O}
_{L}\left( -2\right) \right) .
\end{equation*}
so that $\left( \ref{F1}\right) $ continues to hold for this case.

\subsection{Maps to the $d$-th osculation space}

Let $Y/S^{\prime }$ be a versal family of $k$-folds in $X/S$ such that
Condition 1 and Condition 2 above hold. For example this is the case if 
\begin{equation*}
d+a\geq \max \left\{ \frac{7n-3k-3}{4},\frac{3n-k+1}{2}\right\} .
\end{equation*}
and $Y/S^{\prime }$ is a versal sub-family of $k$-folds such that $\left( 
\ref{fund}\right) $ holds. We finish the proof of Theorem \ref{MT} by
showing:

\begin{lemma}
\label{inside}Suppose $Y/S^{\prime }$ is a versal family of $k$-folds in $X/S
$ such that Condition 1 and Condition 2 above hold.and 
\begin{equation}
\frac{d\left( d-1\right) }{2}-n\geq \left( 2n-1-d\right) -k=codim\left(
g\left( Y\right) ,\Delta _{d}\right) .  \label{newrel}
\end{equation}
Then 
\begin{equation*}
f\left( Y\right) \subseteq \rho _{d+1}\left( \Delta _{d+1}\right) .
\end{equation*}
That is $f\left( Y_{F}\right) $ lies in the sub-variety cut out by the union
of all lines on $X_{F}$.
\end{lemma}

\begin{proof}
At generic $\left( x,l,F\right) \in g\left( Y\right) $, suppose that 
\begin{equation*}
\left. T_{g\left( Y\right) }\cap T_{\pi _{d}}\right| _{\left( x,l,F\right)
}\nsubseteqq \left. \left( \left( p\circ \pi \right) ^{*}M_{G}^{d}\right)
\right| _{\left( x,l,F\right) }.
\end{equation*}
Then letting $c=codim\left( g\left( Y\right) ,\Delta _{d}\right) $, there is
a vector 
\begin{equation*}
v\in \left. \left( \bigwedge\nolimits^{c}\left( p\circ \pi \right)
^{*}M_{G}^{d}\right) \right| _{\left( x,l,F\right) }
\end{equation*}
such that 
\begin{equation}
\left\langle \left. v\wedge w\wedge z\right| \left. \omega _{\Delta
_{d}}\right| _{\left( x,l,F\right) }\right\rangle \neq 0  \label{dil'}
\end{equation}
where 
\begin{eqnarray*}
w &\in &\bigwedge\nolimits^{N}\left. T_{S}\right| _{F} \\
z &\in &\bigwedge\nolimits^{k}T_{\Delta _{d,F},\left( x,l,F\right) }.
\end{eqnarray*}
From Lemma \ref{LP}iv), $\left( \ref{F1}\right) $, and $\left( \ref{newrel}
\right) $, we have that $\left( \bigwedge\nolimits^{c}\pi
^{*}p^{*}M_{G}^{d}\right) \otimes \omega _{\Delta _{d}}$ is generated by
global sections. So there is an element 
\begin{equation*}
\upsilon \in H^{0}\left( \omega _{\Delta _{d}}\otimes
\bigwedge\nolimits^{c}\left( p\circ \pi _{d}\right) ^{*}M_{G}^{d}\right)
\end{equation*}
such that $\upsilon \mapsto v$ under the map 
\begin{equation*}
H^{0}\left( \omega _{\Delta _{d}}\otimes \bigwedge\nolimits^{c}\left( p\circ
\pi _{d}\right) ^{*}M_{G}^{d}\right) \rightarrow \left. \omega _{\Delta
_{d}}\otimes \bigwedge\nolimits^{c}\left( p\circ \pi _{d}\right)
^{*}M_{G}^{d}\right| _{\left( x,l,F\right) }.
\end{equation*}
Let $\alpha $ denote the image of $\upsilon $ under the map 
\begin{equation*}
H^{0}\left( \omega _{\Delta _{d}}\otimes \bigwedge\nolimits^{c}\left( p\circ
\pi _{d}\right) ^{*}M_{G}^{d}\right) \rightarrow H^{0}\left(
\bigwedge\nolimits^{c}T_{\Delta _{d}}\otimes \omega _{\Delta _{d}}\right) 
\overset{\left\langle \left. w\right| \ \right\rangle }{\longrightarrow }
H^{0}\left( \Omega _{\Delta _{d,F}}^{k}\right) \rightarrow H^{0}\left(
\omega _{Y_{F}}^{k}\right) .
\end{equation*}
The hypothesis of the lemma implies that, if $\left( x,l,F\right) $ is
generic in $g\left( Y\right) $%
\begin{equation*}
\left\langle \left. z\right| \alpha _{\left( x,l,F\right) }\right\rangle =0
\end{equation*}
contradicting $\left( \ref{dil'}\right) .$ So 
\begin{equation}
\left. T_{g\left( Y\right) }\right| _{\left( x,l,F\right) }\cap \left.
T_{\pi _{d}}\right| _{\left( x,l,F\right) }\subseteq \left. \left( \left(
p\circ \pi _{d}\right) ^{*}M_{G}^{d}\right) \right| _{\left( x,l,F\right) }.
\label{F4}
\end{equation}

We claim that 
\begin{equation*}
F\in \left( s\circ \rho _{d}\right) _{*}\left( \left. T_{g\left( Y\right)
}\right| _{\left( x,l,F\right) }\cap \left. T_{\pi _{d}}\right| _{\left(
x,l,F\right) }\right) \subseteq \left. M_{G}^{d}\right| _{l}
\end{equation*}
which immediately gives 
\begin{equation*}
\left( l,x,F\right) \in \Delta _{d+1}.
\end{equation*}
To see this notice that $GL\left( n+1\right) $ acts on the diagram 
\begin{equation*}
\begin{array}{ccccc}
Y &  & \overset{f}{\longrightarrow } &  & X \\ 
\downarrow g &  &  & \nearrow _{\rho } & \downarrow \left( t,s\right) \\ 
\Delta _{d} & \rightarrow & L\times _{\Bbb{P}^{n}}S & \rightarrow & \Bbb{P}
^{n}\times S
\end{array}
\end{equation*}
so that the stabilizer of $\Bbb{C}\cdot x$ acts on the map 
\begin{equation*}
\left( s\circ f\right) :\left( t\circ f\right) ^{-1}\left( x\right)
\rightarrow S.
\end{equation*}
This in turn implies the following containment at $\left( l,x,F\right) :$%
\begin{eqnarray*}
\left( F+\sum\nolimits_{i}\left( \frac{\partial F}{\partial X_{i}}%
M_{x}^{1}\right) \right) &\subseteq &\left( s\circ f\right) _{*}\left(
T_{\left( t\circ f\right) ^{-1}\left( x\right) }\right) \\
&=&\left( s\circ \rho _{d}\circ g\right) _{*}\left( T_{\left( t\circ
f\right) ^{-1}\left( x\right) }\right) \\
&=&\left( s\circ \rho _{d}\right) _{*}\left( T_{g\left( \left( q\circ \pi
_{d}\circ g\right) ^{-1}\left( x\right) \right) }\right) \\
&=&\left( s\circ \rho _{d}\right) _{*}\left( T_{g\left( Y\right) }\cap
T_{q\circ \pi _{d}}\right)
\end{eqnarray*}
where $M_{x}^{1}$ denotes the linear forms vanishing on $x$. On the other
hand, the containment 
\begin{equation*}
\left( \rho _{d}\right) _{*}\left( T_{g\left( Y\right) }\cap T_{\pi
_{d}}\right) \subseteq \left( \rho _{d}\right) _{*}\left( T_{g\left(
Y\right) }\cap T_{q\circ \pi _{d}}\right)
\end{equation*}
is actually an equality at generic $\left( l,x,F\right) \in g\left( Y\right) 
$ since the composition 
\begin{equation*}
Y\rightarrow \Delta _{d}\subseteq L\times _{\Bbb{P}^{n}}X\subseteq L\times
X\rightarrow X
\end{equation*}
is an immersion there. So 
\begin{equation*}
F\in \left( s\circ \rho _{d}\right) _{*}\left( T_{g\left( Y\right) }\cap
T_{\pi _{d}}\right) =\left. M_{G}^{d}\right| _{l}.
\end{equation*}
But 
\begin{equation*}
F\in \left. M_{G}^{d}\right| _{l}
\end{equation*}
implies that 
\begin{equation*}
\left. F\right| _{l}=0.
\end{equation*}
\end{proof}

Corollary \ref{MC} is then obtained as follows.

\begin{corollary}
Suppose $\left( n-1\right) -k=2$, 
\begin{equation*}
\frac{d\left( d-1\right) }{2}-n\geq n+2-d=codim\left( g\left( Y\right) ,%
\Delta _{d}\right) ,
\end{equation*}

and 
\begin{equation*}
f\left( Y\right) \nsubseteqq \rho _{d+1}\left( \Delta _{d+1}\right) .
\end{equation*}
Then, for $r$ as in $\left( \ref{rdef}\right) $, 
\begin{equation*}
h^{0}\left( \omega _{Y_{F}}\left( n+2-d\right) \right) \geq \left( 
\begin{array}{c}
n+1 \\ 
k+1
\end{array}
\right) -r.
\end{equation*}
\end{corollary}

\begin{proof}
Since 
\begin{equation*}
s_{T_{t\circ f}}\leq \left( n-1\right) -k=2
\end{equation*}
the only possibilities are 
\begin{equation*}
s_{T_{t\circ f}}=s_{T_{t\circ f}}^{\prime }=1
\end{equation*}
and 
\begin{equation*}
s_{T_{t\circ f}}=2,\ s_{T_{t\circ f}}^{\prime }=0.
\end{equation*}
In the latter case, Theorem \ref{osc} and Lemma \ref{inside} imply that 
\begin{equation*}
f\left( Y\right) \subseteq \rho _{d+1}\left( \Delta _{d+1}\right) .
\end{equation*}
So $s_{T_{t\circ f}}=1$. Now apply Lemma \ref{s=1}.
\end{proof}

\end{document}